\documentclass[12pt,a4paper]{article}
\usepackage{theorem,enumerate}
\usepackage{amsmath,latexsym,amssymb,amsfonts}
\usepackage{eucal}
\usepackage{color}
\usepackage{comment}
\usepackage{eufrak}
\newcounter{Scounter}
\theorembodyfont{\normalfont\slshape}
\newtheorem{thm}{Theorem}[section]

\newtheorem{Thm}{Theorem}

\newtheorem{prop}[thm]{Proposition}
\newtheorem{lem}[thm]{Lemma}

\newtheorem{Q}[thm]{Question} 

\newtheorem{claim}{Claim}[section] 

\newtheorem{fact}[claim]{Fact} 
\newtheorem{con}{Conjecture}

\numberwithin{equation}{section}

\newcommand{\proof}{\medbreak\noindent\textit{Proof.}\quad}

\newcommand{\qed}{{$\quad\square$\vs{3.6}}}

\newcommand{\vs}[1]{\vspace*{#1 mm}}

\def\A{{ \mathcal{A}}}

\def\C{{ \mathcal{C}}}
\def\D{{ \mathcal{D}}}

\def\H{{ \mathcal{H}}}

\def\L{{ \mathcal{L}}}
\def\P{{ \mathcal{P}}}
\def\Q{{ \mathcal{Q}}}

\def\S{{ \mathcal{S}}}

\bfseries\normalfont

\title{The existence of a path-factor without small odd paths}
\author{
Yoshimi Egawa$^1$ \and\
Michitaka Furuya$^1$\footnote{\texttt{e-mail:michitaka.furuya@gmail.com}} \vs{5}\\
$^1$\textsl{Department of Mathematical Information Science,} \\
\textsl{Tokyo University of Science,}\\
\textsl{1-3 Kagurazaka, Shinjuku-ku, Tokyo 162-8601, Japan }\\
}

\date{}

\begin{document}

\maketitle

\begin{abstract}
In this paper, we show that if a graph $G$ satisfies $c_{1}(G-X)+\frac{2}{3}c_{3}(G-X)\leq \frac{4}{3}|X|+\frac{1}{3}$ for all $X\subseteq V(G)$, then $G$ has a $\{P_{2},P_{5}\}$-factor, where $c_{i}(G-X)$ is the number of components $C$ of $G-X$ with $|V(C)|=i$.
\end{abstract}



%%%%%%%%%%%%%%%%%%%%%%%%%%%%%%%%%%%%%%%%%%%%%%%%%%%%%%%%%%%%%%%%%%%%%%%%%%%%%%%%%%%%%%%%%%%%%%%%%%%%%%%%%%%%%%%%%%%%%%%%
%%%%%%%%%%%%%%%%%%%%%%%%%%%%%%%%%%%%%%%%%%%%%%%%%%%%%%%%%%%%%%%%%%%%%%%%%%%%%%%%%%%%%%%%%%%%%%%%%%%%%%%%%%%%%%%%%%%%%%%%
\section{Introduction}\label{sec1}
%%%%%%%%%%%%%%%%%%%%%%%%%%%%%%%%%%%%%%%%%%%%%%%%%%%%%%%%%%%%%%%%%%%%%%%%%%%%%%%%%%%%%%%%%%%%%%%%%%%%%%%%%%%%%%%%%%%%%%%%
%%%%%%%%%%%%%%%%%%%%%%%%%%%%%%%%%%%%%%%%%%%%%%%%%%%%%%%%%%%%%%%%%%%%%%%%%%%%%%%%%%%%%%%%%%%%%%%%%%%%%%%%%%%%%%%%%%%%%%%%

In this paper, all graphs are finite and simple.
Let $G$ be a graph.
We let $V(G)$ and $E(G)$ denote the vertex set and the edge set of $G$, respectively.
For $u\in V(G)$, we let $N_{G}(u)$ and $d_{G}(u)$ denote the {\it neighborhood} and the {\it degree} of $u$, respectively.
For $U\subseteq V(G)$, we let $N_{G}(U)=(\bigcup _{u\in U}N_{G}(u))-U$.
For disjoint sets $X,Y\subseteq V(G)$, we let $E_{G}(X,Y)$ denote the set of edges of $G$ joining a vertex in $X$ and a vertex in $Y$.
For $X\subseteq V(G)$, we let $G[X]$ denote the subgraph of $G$ induced by $X$.
%For a graph $H$ and an integer $s\geq 2$, we let $sH$ denote the disjoint union of $s$ copies of $H$.
For two graphs $H_{1}$ and $H_{2}$, we let $H_{1}+H_{2}$ denote the {\it join} of $H_{1}$ and $H_{2}$.
Let $P_{n}$ denote the {\it path} of order $n$.
For terms and symbols not defined here, we refer the reader to \cite{D}.

For a set $\H$ of connected graphs, a spanning subgraph $F$ of a graph is called an {\it $\H$-factor} if each component of $F$ is isomorphic to a graph in $\H$.
A {\it path-factor} of a graph is a spanning subgraph whose components are paths of order at least $2$.
Since every path of order at least $2$ can be partitioned into paths of orders $2$ and $3$, a graph has a path-factor if and only if it has a $\{P_{2},P_{3}\}$-factor.
Akiyama, Avis and Era~\cite{AAE} gave a necessary and sufficient condition for the existence of a path-factor (here $i(G)$ denotes the number of isolated vertices of a graph $G$).

\begin{Thm}[Akiyama, Avis and Era~\cite{AAE}]%%%%%%%%%%%%%%%%%%%%%%%%%%%%%%%%%%%%%%%%%%%%%%%%%%%%%%%%%%%%%%%%%%%%%%%%%%%
\label{ThmA}
A graph $G$ has a $\{P_{2},P_{3}\}$-factor if and only if $i(G-X)\leq 2|X|$ for all $X\subseteq V(G)$.
\end{Thm}
%%%%%%%%%%%%%%%%%%%%%%%%%%%%%%%%%%%%%%%%%%%%%%%%%%%%%%%%%%%%%%%%%%%%%%%%%%%%%%%%%%%%%%%%%%%%%%%%%%%%%%%%%%%%%%%%%%%%%%%%

Now we consider a path-factor with additional conditions.
For example, one may require a path-factor to consist of components of large order.
Concerning such a problem, Kaneko~\cite{K} gave a necessary and sufficient condition for the existence of a path-factor whose components have order at least $3$.
On the other hand, for $k\geq 4$, it is not known that whether the existence problem of a path-factor whose components have order at least $k$ is polynomially solvable or not, though some results about such a factor have been obtained (see, for example, Kano, Lee and Suzuki~\cite{KLS} and Kawarabayashi, Matsuda, Oda and Ota~\cite{KMOO}).

In this paper, we study a different type of path-factor problem.
Specifically, we focus on the existence of a $\{P_{2},P_{2k+1}\}$-factor ($k\geq 2$).

There are two motivations to study such factors.
One of the motivations is related the notion of a hypomatchable graph.
A graph $H$ is {\it hypomatchable} if $H-x$ has a perfect matching for every $x\in V(H)$.
A graph is a {\it propeller} if it is obtained from a hypomatchable graph $H$ by adding new vertices $a,b$ together with edge $ab$, and joining $a$ to some vertices of $H$.
Loebal and Poljak~\cite{LP} proved the following theorem.

\begin{Thm}[Loebal and Poljak~\cite{LP}]%%%%%%%%%%%%%%%%%%%%%%%%%%%%%%%%%%%%%%%%%%%%%%%%%%%%%%%%%%%%%%%%%%%%%%%%%%%%%%%%
\label{ThmB}
Let $H$ be a connected graph.
If either $H$ has a perfect matching, or $H$ is hypomatchable, or $H$ is a propeller, then the existence problem of a $\{P_{2},H\}$-factor is polynomially solvable.
The problem is {\bf NP}-complete for all other graphs $H$.
\end{Thm}
%%%%%%%%%%%%%%%%%%%%%%%%%%%%%%%%%%%%%%%%%%%%%%%%%%%%%%%%%%%%%%%%%%%%%%%%%%%%%%%%%%%%%%%%%%%%%%%%%%%%%%%%%%%%%%%%%%%%%%%%

In particular, for $k\geq 2$, the existence problem of a $\{P_{2},P_{2k+1}\}$-factor is {\bf NP}-complete.
Because of this fact, existence problems concerning $\{P_{2},P_{2k+1}\}$-factors seem to have unjustly been ignored.
However, in general, the fact that a problem is {\bf NP}-complete in terms of algorithm does not mean that one cannot obtain a theoretical result concerning the problem.
From this viewpoint, in this paper, we prove a theorem on the existence of a $\{P_{2},P_{5}\}$-factor which, we hope, will serve as an initial attempt to develop the theory of $\{P_{2},P_{2k+1}\}$-factors.

The other motivation is the fact that a $\{P_{2},P_{2k+1}\}$-factor is a useful tool for finding large matchings.
It is easy to see that if a graph $G$ has a $\{P_{2},P_{2k+1}\}$-factor, then $G$ has a matching $M$ with $|M|\geq \frac{k}{2k+1}|V(G)|$.
Thus the existence of a $\{P_{2},P_{2k+1}\}$-factor helps to find large matchings.
%, and hence we cannot expect a manageable necessary and sufficient condition for it.
%Instead of that, we seek a better sufficient condition for it.

In order to state our theorem, we need some more definitions.
For a graph $H$, we let $\C(H)$ be the set of components of $H$, and for $i\geq 1$, let $\C _{i}(H)=\{C\in \C(H)\mid |V(C)|=i\}$ and $c_{i}(H)=|\C _{i}(H)|$.
Note that $c_{1}(H)$ is the number of isolated vertices of $H$ (i.e., $c_{1}(H)=i(H)$).
If a graph $G$ has a $\{P_{2},P_{5}\}$-factor, then $c_{1}(G-X)+\frac{1}{2}c_{3}(G-X)\leq \frac{3}{2}|X|$ for all $X\subseteq V(G)$ (see Section~\ref{sec-nec}).
Thus if a condition concerning $c_{1}(G-X)$ and $c_{3}(G-X)$ for $X\subseteq V(G)$ assures us the existence of a $\{P_{2},P_{5}\}$-factor, then it will make a useful sufficient condition.

The main purpose of this paper is to prove the following theorem.

\begin{thm}%%%%%%%%%%%%%%%%%%%%%%%%%%%%%%%%%%%%%%%%%%%%%%%%%%%%%%%%%%%%%%%%%%%%%%%%%%%%%%%%%%%%%%%%%%%%%%%%%%%%%%%%%%%%%
\label{thm1}
Let $G$ be a graph.
If $c_{1}(G-X)+\frac{2}{3}c_{3}(G-X)\leq \frac{4}{3}|X|+\frac{1}{3}$ for all $X\subseteq V(G)$, then $G$ has a $\{P_{2},P_{5}\}$-factor.
\end{thm}
%%%%%%%%%%%%%%%%%%%%%%%%%%%%%%%%%%%%%%%%%%%%%%%%%%%%%%%%%%%%%%%%%%%%%%%%%%%%%%%%%%%%%%%%%%%%%%%%%%%%%%%%%%%%%%%%%%%%%%%%

We prove Theorem~\ref{thm1} in Sections~\ref{sec2} and \ref{sec3}.
In Subsection~\ref{sec4.1}, we show that the bound $\frac{4}{3}|X|+\frac{1}{3}$ in Theorem~\ref{thm1} is best possible.

In our proof of Theorem~\ref{thm1}, we make use of the following fact.

\begin{fact}%%%%%%%%%%%%%%%%%%%%%%%%%%%%%%%%%%%%%%%%%%%%%%%%%%%%%%%%%%%%%%%%%%%%%%%%%%%%%%%%%%%%%%%%%%%%%%%%%%%%%%%%%%%%
\label{fact1}
Let $G$ be a graph.
Then $G$ has a $\{P_{2},P_{5}\}$-factor if and only if $G$ has a path-factor $F$ with $\C_{3}(F)=\emptyset $.
\end{fact}
%%%%%%%%%%%%%%%%%%%%%%%%%%%%%%%%%%%%%%%%%%%%%%%%%%%%%%%%%%%%%%%%%%%%%%%%%%%%%%%%%%%%%%%%%%%%%%%%%%%%%%%%%%%%%%%%%%%%%%%%

We conclude this section with a conjecture concerning $\{P_{2},P_{2k+1}\}$-factors with $k\geq 3$.
By Theorems~\ref{ThmA} and \ref{thm1}, for $k\in \{1,2\}$, there exists a constant $a_{k}>1$ such that the condition $\sum _{0\leq i\leq k-1}c_{2i+1}(G-X)\leq a_{k}|X|~(X\subseteq V(G))$ assures us the existence of a $\{P_{2},P_{2k+1}\}$-factor (one can take $a_{1}=2$ and $a_{2}=\frac{4}{3}$).
Thus one may expect that there exists a similar constant $a_{k}>1$ for $k\geq 3$.
However, when we consider the case where $k\geq 3$ with $k\equiv 0~(\mbox{mod }3)$, the situation changes drastically; that is, there exist infinitely many graphs $G$ having no $\{P_{2},P_{2k+1}\}$-factor such that $\sum _{0\leq i\leq k-1}c_{2i+1}(G-X)\leq \frac{4k+6}{8k+3}|X|+\frac{2k+3}{8k+3}$ for all $X\subseteq V(G)$ (see Subsection~\ref{sec4.2}).
Thus we pose the following conjecture.

\begin{con}%%%%%%%%%%%%%%%%%%%%%%%%%%%%%%%%%%%%%%%%%%%%%%%%%%%%%%%%%%%%%%%%%%%%%%%%%%%%%%%%%%%%%%%%%%%%%%%%%%%%%%%%%%%%%
\label{con1}
Let $k\geq 3$, and let $G$ be a graph.
If $\sum _{0\leq i\leq k-1}c_{2i+1}(G-X)\leq \frac{4k+6}{8k+3}|X|$ for all $X\subseteq V(G)$, then $G$ has a $\{P_{2},P_{2k+1}\}$-factor.
\end{con}
%%%%%%%%%%%%%%%%%%%%%%%%%%%%%%%%%%%%%%%%%%%%%%%%%%%%%%%%%%%%%%%%%%%%%%%%%%%%%%%%%%%%%%%%%%%%%%%%%%%%%%%%%%%%%%%%%%%%%%%%

%%%%%%%%%%%%%%%%%%%%%%%%%%%%%%%%%%%%%%%%%%%%%%%%%%%%%%%%%%%%%%%%%%%%%%%%%%%%%%%%%%%%%%%%%%%%%%%%%%%%%%%%%%%%%%%%%%%%%%%%
%%%%%%%%%%%%%%%%%%%%%%%%%%%%%%%%%%%%%%%%%%%%%%%%%%%%%%%%%%%%%%%%%%%%%%%%%%%%%%%%%%%%%%%%%%%%%%%%%%%%%%%%%%%%%%%%%%%%%%%%
\section{A necessary condition for a $\{P_{2},P_{5}\}$-factor}\label{sec-nec}
%%%%%%%%%%%%%%%%%%%%%%%%%%%%%%%%%%%%%%%%%%%%%%%%%%%%%%%%%%%%%%%%%%%%%%%%%%%%%%%%%%%%%%%%%%%%%%%%%%%%%%%%%%%%%%%%%%%%%%%%
%%%%%%%%%%%%%%%%%%%%%%%%%%%%%%%%%%%%%%%%%%%%%%%%%%%%%%%%%%%%%%%%%%%%%%%%%%%%%%%%%%%%%%%%%%%%%%%%%%%%%%%%%%%%%%%%%%%%%%%%

In this section, we give a necessary condition for the existence of a $\{P_{2},P_{5}\}$-factor in terms of invariants $c_{1}$ and $c_{3}$.
We show the following proposition.

\begin{prop}%%%%%%%%%%%%%%%%%%%%%%%%%%%%%%%%%%%%%%%%%%%%%%%%%%%%%%%%%%%%%%%%%%%%%%%%%%%%%%%%%%%%%%%%%%%%%%%%%%%%%%%%%%%%
\label{prop1.1}
If a graph $G$ has a $\{P_{2},P_{5}\}$-factor, then $c_{1}(G-X)+\frac{1}{2}c_{3}(G-X)\leq \frac{3}{2}|X|$ for all $X\subseteq V(G)$.
\end{prop}
%%%%%%%%%%%%%%%%%%%%%%%%%%%%%%%%%%%%%%%%%%%%%%%%%%%%%%%%%%%%%%%%%%%%%%%%%%%%%%%%%%%%%%%%%%%%%%%%%%%%%%%%%%%%%%%%%%%%%%%%
\proof
Let $F$ be a $\{P_{2},P_{5}\}$-factor of $G$, and let $X\subseteq V(G)$.
Then we can verify that
\begin{align}
c_{1}(P-X)+\frac{1}{2}c_{3}(P-X)\leq \frac{3}{2}|V(P)\cap X| \mbox{ for every } P\in \C(F).\label{sec1-1}
\end{align}
Since every component $C$ of $G-X$ with $|V(C)|=1$ belongs to $\bigcup _{P\in \C(F)}\C_{1}(P-X)$, we have
\begin{align}
|&\C_{1}(G-X)|=\sum _{P\in \C(F)}|\C_{1}(P-X)|-\left|\left(\bigcup _{P\in \C(F)}\C_{1}(P-X)\right)-\C_{1}(G-X)\right|.\label{sec1-2}
\end{align}
Furthermore,
\begin{align}
|&\C_{3}(G-X)|\leq \sum _{P\in \C(F)}|\C_{3}(P-X)|+\left|\C_{3}(G-X)-\left(\bigcup _{P\in \C(F)}\C_{3}(P-X)\right)\right|.\label{sec1-3}
\end{align}

Let $C$ be a component of $G-X$ with $|V(C)|=3$ which does not belong to $\bigcup _{P\in \C(F)}\C_{3}(P-X)$.
Then $C$ intersects with at least two components of $F-X$.
Since $|V(C)|=3$, $C$ contains a component of $P-X$ of order $1$ for some $P\in \C(F)$.
Since $C$ is arbitrary, this implies that
\begin{align}
\left|\C_{3}(G-X)-\left(\bigcup _{P\in \C(F)}\C_{3}(P-X)\right)\right|\leq \left|\left(\bigcup _{P\in \C(F)}\C_{1}(P-X)\right)-\C_{1}(G-X)\right|.\label{sec1-4}
\end{align}
By (\ref{sec1-1})--(\ref{sec1-4}),
\begin{align}
c_{1}(G-X)&+\frac{1}{2}c_{3}(G-X)\nonumber \\
&\leq \left(\sum _{P\in \C(F)}|\C_{1}(P-X)|-\left|\left(\bigcup _{P\in \C(F)}\C_{1}(P-X)\right)-\C_{1}(G-X)\right|\right)\nonumber \\
&\quad +\frac{1}{2}\left(\sum _{P\in \C(F)}|\C_{3}(P-X)|+\left|\C_{3}(G-X)-\left(\bigcup _{P\in \C(F)}\C_{3}(P-X)\right)\right|\right)\nonumber \\
&\leq \left(\sum _{P\in \C(F)}|\C_{1}(P-X)|-\left|\left(\bigcup _{P\in \C(F)}\C_{1}(P-X)\right)-\C_{1}(G-X)\right|\right)\nonumber \\
&\quad +\frac{1}{2}\left(\sum _{P\in \C(F)}|\C_{3}(P-X)|+\left|\left(\bigcup _{P\in \C(F)}\C_{1}(P-X)\right)-\C_{1}(G-X)\right|\right)\nonumber \\
&\leq \sum _{P\in \C(F)}|\C_{1}(P-X)|+\frac{1}{2}\sum _{P\in \C(F)}|\C_{3}(P-X)|\nonumber \\
&= \sum _{P\in \C(F)}\left(c_{1}(P-X)+\frac{1}{2}c_{3}(P-X)\right)\nonumber \\
&\leq \frac{3}{2}\sum _{P\in \C(F)}|V(P)\cap X|\nonumber \\
&= \frac{3}{2}|X|.\nonumber
\end{align}
Thus we get the desired conclusion.
\qed

%%%%%%%%%%%%%%%%%%%%%%%%%%%%%%%%%%%%%%%%%%%%%%%%%%%%%%%%%%%%%%%%%%%%%%%%%%%%%%%%%%%%%%%%%%%%%%%%%%%%%%%%%%%%%%%%%%%%%%%%
%%%%%%%%%%%%%%%%%%%%%%%%%%%%%%%%%%%%%%%%%%%%%%%%%%%%%%%%%%%%%%%%%%%%%%%%%%%%%%%%%%%%%%%%%%%%%%%%%%%%%%%%%%%%%%%%%%%%%%%%
\section{A path-factor in bipartite graph}\label{sec2}
%%%%%%%%%%%%%%%%%%%%%%%%%%%%%%%%%%%%%%%%%%%%%%%%%%%%%%%%%%%%%%%%%%%%%%%%%%%%%%%%%%%%%%%%%%%%%%%%%%%%%%%%%%%%%%%%%%%%%%%%
%%%%%%%%%%%%%%%%%%%%%%%%%%%%%%%%%%%%%%%%%%%%%%%%%%%%%%%%%%%%%%%%%%%%%%%%%%%%%%%%%%%%%%%%%%%%%%%%%%%%%%%%%%%%%%%%%%%%%%%%

Let $G$ be a bipartite graph with bipartition $(S,T)$.
A subgraph $F$ of $G$ is {\it $S$-central} if $S\subseteq V(F)$ and $|V(A)\cap T|\geq |V(A)\cap S|$ for every $A\in \C(F)$.

In this section, we focus on the existence of a special path-factor in bipartite graphs, and show the following theorem, which will be used in our proof of Theorem~\ref{thm1}.

\begin{thm}%%%%%%%%%%%%%%%%%%%%%%%%%%%%%%%%%%%%%%%%%%%%%%%%%%%%%%%%%%%%%%%%%%%%%%%%%%%%%%%%%%%%%%%%%%%%%%%%%%%%%%%%%%%%%
\label{thm2}
Let $S$, $T_{1}$ and $T_{2}$ be disjoint sets with $1\leq |S|\leq |T_{1}|+|T_{2}|$ and $|T_{1}|+\frac{2}{3}|T_{2}|\leq \frac{4}{3}|S|+\frac{1}{3}$, and set $T=T_{1}\cup T_{2}$.
Let $G$ be a bipartite graph with bipartition $(S,T)$ satisfying the property that for every $X\subseteq V(G)$, we have either $|N_{G}(X)\cap T_{1}|+\frac{2}{3}|N_{G}(X)\cap T_{2}|\geq \frac{4}{3}|X|$ or $N_{G}(X)=T$.
Then $G$ has an $S$-central path-factor $F$ such that $V(A)\cap T_{2}\not=\emptyset $ for every $A\in \C_{3}(F)$.
\end{thm}
%%%%%%%%%%%%%%%%%%%%%%%%%%%%%%%%%%%%%%%%%%%%%%%%%%%%%%%%%%%%%%%%%%%%%%%%%%%%%%%%%%%%%%%%%%%%%%%%%%%%%%%%%%%%%%%%%%%%%%%%

Before proving the theorem, we prove a lemma.

\begin{lem}%%%%%%%%%%%%%%%%%%%%%%%%%%%%%%%%%%%%%%%%%%%%%%%%%%%%%%%%%%%%%%%%%%%%%%%%%%%%%%%%%%%%%%%%%%%%%%%%%%%%%%%%%%%%%
\label{lem2.1}
Let $S$, $T_{1}$, $T_{2}$, $T$ and $G$ be as in Theorem~\ref{thm2}.
Then $G$ has an $S$-central path-factor.
\end{lem}
%%%%%%%%%%%%%%%%%%%%%%%%%%%%%%%%%%%%%%%%%%%%%%%%%%%%%%%%%%%%%%%%%%%%%%%%%%%%%%%%%%%%%%%%%%%%%%%%%%%%%%%%%%%%%%%%%%%%%%%%
\proof
Let $X\subseteq S$.
If $|N_{G}(X)\cap T_{1}|+\frac{2}{3}|N_{G}(X)\cap T_{2}|\geq \frac{4}{3}|X|$, then $|N_{G}(X)|\geq |N_{G}(X)\cap T_{1}|+\frac{2}{3}|N_{G}(X)\cap T_{2}|\geq \frac{4}{3}|X|\geq |X|$; if $N_{G}(X)=T_{1}\cup T_{2}$, then $|N_{G}(X)|=|T_{1}|+|T_{2}|\geq |S|\geq |X|$.
In either case, we have $|N_{G}(X)|\geq |X|$.
Since $X$ is arbitrary, $G$ has a matching covering $S$ by Hall's marriage theorem.
In particular, $G$ has an $S$-central subgraph $F$ such that every component of $F$ is a path of order at least $2$.
Choose $F$ so that $|V(F)|$ is as large as possible.

Suppose that $V(G)-V(F)\not=\emptyset $.
Note that $V(G)-V(F)\subseteq T$.
Now we define the set $\A$ of components of $F$ as follows:
Let $\A_{1}$ be the set of components $A$ of $F$ with $E_{G}(V(A)\cap S,V(G)-V(F))\not=\emptyset $.
For each $i\geq 2$, let $\A_{i}$ be the set of components $A$ of $F$ with $A\not\in \bigcup _{1\leq j\leq i-1}\A_{j}$ and $E_{G}(V(A)\cap S,\bigcup _{A'\in \A_{i-1}}(V(A')\cap T))\not=\emptyset $.
Let $\A=\bigcup _{i\geq 1}\A_{i}$.

\begin{claim}%%%%%%%%%%%%%%%%%%%%%%%%%%%%%%%%%%%%%%%%%%%%%%%%%%%%%%%%%%%%%%%%%%%%%%%%%%%%%%%%%%%%%%%%%%%%%%%%%%%%%%%%%%%
\label{cl2.1.1}
Every path belonging to $\A$ is isomorphic to $P_{3}$.
\end{claim}
%%%%%%%%%%%%%%%%%%%%%%%%%%%%%%%%%%%%%%%%%%%%%%%%%%%%%%%%%%%%%%%%%%%%%%%%%%%%%%%%%%%%%%%%%%%%%%%%%%%%%%%%%%%%%%%%%%%%%%%%
\proof
Suppose that $\A$ contains a path which is not isomorphic to $P_{3}$.
Let $i$ be the minimum integer such that $\A_{i}$ contains a path $A_{i}=v^{(i)}_{1}\cdots v^{(i)}_{l}$ with $A_{i}\not\simeq P_{3}$.
By the minimality of $i$, every path belonging to $\bigcup _{1\leq j\leq i-1}\A_{j}$ is isomorphic to $P_{3}$.
Hence by the definition of $\A_{j}$, there exists a vertex $v^{(0)}_{1}\in V(G)-V(F)$ and there exist paths $A_{j}=v^{(j)}_{1}v^{(j)}_{2}v^{(j)}_{3}\in \A_{j}~(1\leq j\leq i-1)$ such that $E_{G}(V(A_{1})\cap S,\{v^{(0)}_{1}\})\not=\emptyset $ and $E_{G}(V(A_{j+1})\cap S,V(A_{j})\cap T)\not=\emptyset $ for every $j~(1\leq j\leq i-1)$.
For each $j~(1\leq j\leq i-1)$, by renumbering the veritices $v^{(j)}_{1},v^{(j)}_{2},v^{(j)}_{3}$ of $A_{j}$ backward (i.e., by tracing the path $v^{(j)}_{1}v^{(j)}_{2}v^{(j)}_{3}$ backward and numbering the vertices accordingly) if necessary, we may assume that $E_{G}(V(A_{j+1})\cap S,\{v^{(j)}_{1}\})\not=\emptyset $.
Let $m$ be an index such that $v^{(i)}_{m}v^{(i-1)}_{1}\in E(G)$.
Note that $l\geq 2$ and $l\not=3$.
Thus by renumbering the vertices $v^{(i)}_{1},\ldots ,v^{(i)}_{l}$ of $A_{i}$ backward if necessary, we may assume that $m\not=2$ if $l$ is odd, and $m$ is odd if $l$ is even.
Let $B_{j}=v^{(j-1)}_{1}v^{(j)}_{2}v^{(j)}_{3}~(1\leq j\leq i-1)$, $B_{i}=v^{(i-1)}_{1}v^{(i)}_{m}v^{(i)}_{m+1}\cdots v^{(i)}_{l}$ and $B_{i+1}=v^{(i)}_{1}\cdots v^{(i)}_{m-1}$ (note that $B_{i+1}=\emptyset $ if and only if $l$ is even and $m=1$).
Then $|V(B_{j})\cap T|\geq |V(B_{j})\cap S|$ for every $j~(1\leq j\leq i+1)$.
Therefore $F'=(F-(\bigcup _{1\leq j\leq i}V(A_{j})))\cup (\bigcup _{1\leq j\leq i+1}B_{j})$ is an $S$-central subgraph of $G$ such that $V(F')=V(F)\cup \{v^{(0)}_{1}\}$ and every component of $F'$ is a path of order at least $2$, which contradicts the maximality of $F$.
\qed

We continue with the proof of the lemma.
Let $X_{0}=(\bigcup _{A\in \A}V(A))\cap S$ and $Y_{0}=((\bigcup _{A\in \A}V(A))\cap T)\cup (V(G)-V(F))$.
Since $V(G)-V(F)\not=\emptyset $ and $\A\subseteq \C_{3}(F)$ by Claim~\ref{cl2.1.1}, we have
\begin{align}
|Y_{0}\cap T_{1}|+\frac{2}{3}|Y_{0}\cap T_{2}|\geq \frac{2}{3}|Y_{0}|\geq \frac{2}{3}(2|X_{0}|+1).\label{lem2.1-1}
\end{align}
By the definition of $\A$, $N_{G}(S-X_{0})\cap Y_{0}=\emptyset $.
In particular, $N_{G}(S-X_{0})\not=T$, and hence $|N_{G}(S-X_{0})\cap T_{1}|+\frac{2}{3}|N_{G}(S-X_{0})\cap T_{2}|\geq \frac{4}{3}|S-X_{0}|$.
This together with (\ref{lem2.1-1}) implies that
\begin{eqnarray*}
|T_{1}|+\frac{2}{3}|T_{2}| &\geq & (|Y_{0}\cap T_{1}|+|N_{G}(S-X_{0})\cap T_{1}|)+\frac{2}{3}(|Y_{0}\cap T_{2}|+|N_{G}(S-X_{0})\cap T_{2}|)\\
&\geq & \frac{2}{3}(2|X_{0}|+1)+\frac{4}{3}|S-X_{0}|\\
&=& \frac{4}{3}|S|+\frac{2}{3},
\end{eqnarray*}
which contradicts the assumption that $|T_{1}|+\frac{2}{3}|T_{2}|\leq \frac{4}{3}|S|+\frac{1}{3}$, completing the proof of the lemma.
\qed

We here outline the proof of Theorem~\ref{thm2}.
We choose an $S$-central path-factor $F_{0}$ so that $F_{0}$ will satisfy certain minimality conditions (see the paragraph following the proof of Claim~\ref{cl2.1}).
We then introduce operations which turn $F_{0}$ into a new path-factor (see the paragraphs following Claim~\ref{cl-type2} and Claim~\ref{cl2.2}), and show that the new path-factor contradicts our choice of $F_{0}$.

\medbreak\noindent\textit{Proof of Theorem~\ref{thm2}.}\quad
We start with some definitions.
Let $F$ be an $S$-central path-factor of $G$.
For each integer $i\geq 2$, let $\C^{(1)}_{i}(F)=\{A\in \C_{i}(F)\mid V(A)\cap T_{2}=\emptyset \}$ and $\C^{(2)}_{i}(F)=\C_{i}(F)-\C^{(1)}_{i}(F)$.
If there is no fear of confusion, we simply write $\C_{i}$ and $\C^{(h)}_{i}~(h\in \{1,2\})$ instead of $\C_{i}(F)$ and $\C^{(h)}_{i}(F)$, respectively.

\begin{figure}
\begin{center}
%WinTpicVersion4.26
\unitlength 0.1in
\begin{picture}( 24.0500, 11.0000)( 15.0500,-16.7500)
% CIRCLE 2 0 0 0 Black Black
% 4 1800 1000 1800 1050 1800 1050 1800 1050
% 
\special{sh 1.000}%
\special{ia 1800 1000 50 50  0.0000000  6.2831853}%
\special{pn 8}%
\special{pn 8}%
\special{ar 1800 1000 50 50  0.0000000  6.2831853}%
% CIRCLE 2 0 0 0 Black Black
% 4 1900 1400 1900 1450 1900 1450 1900 1450
% 
\special{sh 1.000}%
\special{ia 1900 1400 50 50  0.0000000  6.2831853}%
\special{pn 8}%
\special{pn 8}%
\special{ar 1900 1400 50 50  0.0000000  6.2831853}%
% CIRCLE 2 0 0 0 Black Black
% 4 2000 1000 2000 1050 2000 1050 2000 1050
% 
\special{sh 1.000}%
\special{ia 2000 1000 50 50  0.0000000  6.2831853}%
\special{pn 8}%
\special{pn 8}%
\special{ar 2000 1000 50 50  0.0000000  6.2831853}%
% CIRCLE 2 0 0 0 Black Black
% 4 2200 1000 2200 1050 2200 1050 2200 1050
% 
\special{sh 1.000}%
\special{ia 2200 1000 50 50  0.0000000  6.2831853}%
\special{pn 8}%
\special{pn 8}%
\special{ar 2200 1000 50 50  0.0000000  6.2831853}%
% CIRCLE 2 0 0 0 Black Black
% 4 2400 1000 2400 1050 2400 1050 2400 1050
% 
\special{sh 1.000}%
\special{ia 2400 1000 50 50  0.0000000  6.2831853}%
\special{pn 8}%
\special{pn 8}%
\special{ar 2400 1000 50 50  0.0000000  6.2831853}%
% CIRCLE 2 0 0 0 Black Black
% 4 2100 1400 2100 1450 2100 1450 2100 1450
% 
\special{sh 1.000}%
\special{ia 2100 1400 50 50  0.0000000  6.2831853}%
\special{pn 8}%
\special{pn 8}%
\special{ar 2100 1400 50 50  0.0000000  6.2831853}%
% CIRCLE 2 0 0 0 Black Black
% 4 2300 1400 2300 1450 2300 1450 2300 1450
% 
\special{sh 1.000}%
\special{ia 2300 1400 50 50  0.0000000  6.2831853}%
\special{pn 8}%
\special{pn 8}%
\special{ar 2300 1400 50 50  0.0000000  6.2831853}%
% CIRCLE 2 0 0 0 Black Black
% 4 2600 1000 2600 1050 2600 1050 2600 1050
% 
\special{sh 1.000}%
\special{ia 2600 1000 50 50  0.0000000  6.2831853}%
\special{pn 8}%
\special{pn 8}%
\special{ar 2600 1000 50 50  0.0000000  6.2831853}%
% CIRCLE 2 0 0 0 Black Black
% 4 2500 1400 2500 1450 2500 1450 2500 1450
% 
\special{sh 1.000}%
\special{ia 2500 1400 50 50  0.0000000  6.2831853}%
\special{pn 8}%
\special{pn 8}%
\special{ar 2500 1400 50 50  0.0000000  6.2831853}%
% CIRCLE 2 0 0 0 Black Black
% 4 3200 1000 3200 1050 3200 1050 3200 1050
% 
\special{sh 1.000}%
\special{ia 3200 1000 50 50  0.0000000  6.2831853}%
\special{pn 8}%
\special{pn 8}%
\special{ar 3200 1000 50 50  0.0000000  6.2831853}%
% CIRCLE 2 0 0 0 Black Black
% 4 3400 1000 3400 1050 3400 1050 3400 1050
% 
\special{sh 1.000}%
\special{ia 3400 1000 50 50  0.0000000  6.2831853}%
\special{pn 8}%
\special{pn 8}%
\special{ar 3400 1000 50 50  0.0000000  6.2831853}%
% CIRCLE 2 0 0 0 Black Black
% 4 3600 1000 3600 1050 3600 1050 3600 1050
% 
\special{sh 1.000}%
\special{ia 3600 1000 50 50  0.0000000  6.2831853}%
\special{pn 8}%
\special{pn 8}%
\special{ar 3600 1000 50 50  0.0000000  6.2831853}%
% CIRCLE 2 0 0 0 Black Black
% 4 3300 1400 3300 1450 3300 1450 3300 1450
% 
\special{sh 1.000}%
\special{ia 3300 1400 50 50  0.0000000  6.2831853}%
\special{pn 8}%
\special{pn 8}%
\special{ar 3300 1400 50 50  0.0000000  6.2831853}%
% CIRCLE 2 0 0 0 Black Black
% 4 3500 1400 3500 1450 3500 1450 3500 1450
% 
\special{sh 1.000}%
\special{ia 3500 1400 50 50  0.0000000  6.2831853}%
\special{pn 8}%
\special{pn 8}%
\special{ar 3500 1400 50 50  0.0000000  6.2831853}%
% CIRCLE 2 0 0 0 Black Black
% 4 3800 1000 3800 1050 3800 1050 3800 1050
% 
\special{sh 1.000}%
\special{ia 3800 1000 50 50  0.0000000  6.2831853}%
\special{pn 8}%
\special{pn 8}%
\special{ar 3800 1000 50 50  0.0000000  6.2831853}%
% CIRCLE 2 0 0 0 Black Black
% 4 3700 1400 3700 1450 3700 1450 3700 1450
% 
\special{sh 1.000}%
\special{ia 3700 1400 50 50  0.0000000  6.2831853}%
\special{pn 8}%
\special{pn 8}%
\special{ar 3700 1400 50 50  0.0000000  6.2831853}%
% LINE 0 0 3 0 Black Black
% 2 1800 1000 1900 1400
% 
\special{pn 20}%
\special{pa 1800 1000}%
\special{pa 1900 1400}%
\special{fp}%
% LINE 0 0 3 0 Black Black
% 2 2000 1000 2100 1400
% 
\special{pn 20}%
\special{pa 2000 1000}%
\special{pa 2100 1400}%
\special{fp}%
% LINE 0 0 3 0 Black Black
% 2 2200 1000 2300 1400
% 
\special{pn 20}%
\special{pa 2200 1000}%
\special{pa 2300 1400}%
\special{fp}%
% LINE 0 0 3 0 Black Black
% 2 2400 1000 2500 1400
% 
\special{pn 20}%
\special{pa 2400 1000}%
\special{pa 2500 1400}%
\special{fp}%
% LINE 0 0 3 0 Black Black
% 2 3200 1000 3300 1400
% 
\special{pn 20}%
\special{pa 3200 1000}%
\special{pa 3300 1400}%
\special{fp}%
% LINE 0 0 3 0 Black Black
% 2 3400 1000 3500 1400
% 
\special{pn 20}%
\special{pa 3400 1000}%
\special{pa 3500 1400}%
\special{fp}%
% LINE 0 0 3 0 Black Black
% 2 3600 1000 3700 1400
% 
\special{pn 20}%
\special{pa 3600 1000}%
\special{pa 3700 1400}%
\special{fp}%
% LINE 0 0 3 0 Black Black
% 2 3600 1000 3500 1400
% 
\special{pn 20}%
\special{pa 3600 1000}%
\special{pa 3500 1400}%
\special{fp}%
% LINE 0 0 3 0 Black Black
% 2 3800 1000 3700 1400
% 
\special{pn 20}%
\special{pa 3800 1000}%
\special{pa 3700 1400}%
\special{fp}%
% LINE 0 0 3 0 Black Black
% 2 3400 1000 3300 1400
% 
\special{pn 20}%
\special{pa 3400 1000}%
\special{pa 3300 1400}%
\special{fp}%
% LINE 0 0 3 0 Black Black
% 2 2600 1000 2500 1400
% 
\special{pn 20}%
\special{pa 2600 1000}%
\special{pa 2500 1400}%
\special{fp}%
% LINE 0 0 3 0 Black Black
% 2 2400 1000 2300 1400
% 
\special{pn 20}%
\special{pa 2400 1000}%
\special{pa 2300 1400}%
\special{fp}%
% LINE 0 0 3 0 Black Black
% 2 2200 1000 2100 1400
% 
\special{pn 20}%
\special{pa 2200 1000}%
\special{pa 2100 1400}%
\special{fp}%
% LINE 0 0 3 0 Black Black
% 2 2000 1000 1900 1400
% 
\special{pn 20}%
\special{pa 2000 1000}%
\special{pa 1900 1400}%
\special{fp}%
% SPLINE 2 0 3 0 Black Black
% 6 2300 1400 2300 1000 2410 850 3210 850 3400 1000 3400 1000
% 
\special{pn 8}%
\special{pa 2300 1400}%
\special{pa 2294 1368}%
\special{pa 2288 1334}%
\special{pa 2284 1300}%
\special{pa 2280 1268}%
\special{pa 2276 1234}%
\special{pa 2272 1202}%
\special{pa 2272 1138}%
\special{pa 2274 1106}%
\special{pa 2280 1076}%
\special{pa 2286 1044}%
\special{pa 2294 1016}%
\special{pa 2306 986}%
\special{pa 2338 930}%
\special{pa 2358 904}%
\special{pa 2380 880}%
\special{pa 2402 858}%
\special{pa 2428 836}%
\special{pa 2454 818}%
\special{pa 2482 800}%
\special{pa 2510 786}%
\special{pa 2540 772}%
\special{pa 2572 762}%
\special{pa 2602 752}%
\special{pa 2636 744}%
\special{pa 2668 738}%
\special{pa 2702 734}%
\special{pa 2736 732}%
\special{pa 2772 730}%
\special{pa 2806 730}%
\special{pa 2840 732}%
\special{pa 2876 736}%
\special{pa 2910 742}%
\special{pa 2946 748}%
\special{pa 3014 764}%
\special{pa 3046 776}%
\special{pa 3080 786}%
\special{pa 3112 800}%
\special{pa 3172 828}%
\special{pa 3202 846}%
\special{pa 3230 862}%
\special{pa 3256 880}%
\special{pa 3308 920}%
\special{pa 3356 960}%
\special{pa 3380 982}%
\special{pa 3400 1000}%
\special{sp}%
% STR 2 0 3 0 Black Black
% 4 1600 1300 1600 1400 5 0 0 0
% $A$
\put(16.0000,-14.0000){\makebox(0,0){$A$}}%
% STR 2 0 3 0 Black Black
% 4 4000 1300 4000 1400 5 0 0 0
% $B$
\put(40.0000,-14.0000){\makebox(0,0){$B$}}%
% STR 2 0 3 0 Black Black
% 4 2300 1640 2300 1740 5 0 0 0
% $\sigma _{F}(AB)$
\put(23.0000,-17.4000){\makebox(0,0){$\sigma _{F}(AB)$}}%
% STR 2 0 3 0 Black Black
% 4 3660 620 3660 720 5 0 0 0
% $\tau _{F}(AB)$
\put(36.6000,-7.2000){\makebox(0,0){$\tau _{F}(AB)$}}%
% STR 2 0 3 0 Black Black
% 4 2600 540 2600 640 5 0 0 0
% $\varphi _{F}(AB)$
\put(26.0000,-6.4000){\makebox(0,0){$\varphi _{F}(AB)$}}%
% VECTOR 2 0 3 0 Black Black
% 2 2300 1630 2300 1480
% 
\special{pn 8}%
\special{pa 2300 1630}%
\special{pa 2300 1480}%
\special{fp}%
\special{sh 1}%
\special{pa 2300 1480}%
\special{pa 2280 1548}%
\special{pa 2300 1534}%
\special{pa 2320 1548}%
\special{pa 2300 1480}%
\special{fp}%
% VECTOR 2 0 3 0 Black Black
% 2 3600 800 3450 940
% 
\special{pn 8}%
\special{pa 3600 800}%
\special{pa 3450 940}%
\special{fp}%
\special{sh 1}%
\special{pa 3450 940}%
\special{pa 3512 910}%
\special{pa 3490 904}%
\special{pa 3486 880}%
\special{pa 3450 940}%
\special{fp}%
\end{picture}%
\caption{Edge $\varphi _{F}(AB)$ and vertices $\sigma _{F}(AB)$ and $\tau _{F}(AB)$}
\label{thm2.1fig1}
\end{center}
\end{figure}
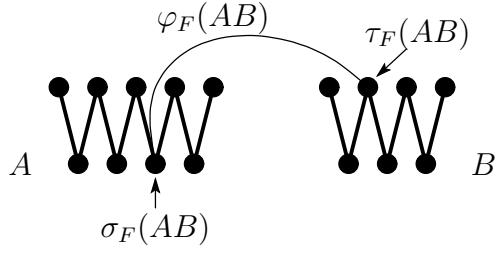

Let $\D_{F}$ be the digraph defined by $V(\D_{F})=\C(F)$ and $E(\D)=\{AB\mid E_{G}(V(A)\cap S,V(B)\cap T)\not=\emptyset \}$.
For each edge $AB\in E(\D_{F})$, we fix an edge $\varphi _{F}(AB)$ in $E_{G}(V(A)\cap S,V(B)\cap T)$, and let $\sigma _{F}(AB)\in V(G)$ be the vertex of $A$ incident with $\varphi _{F}(AB)$ and $\tau _{F}(AB)\in V(G)$ be the vertex of $B$ incident with $\varphi _{F}(AB)$ (see Figure~\ref{thm2.1fig1}).

For a path $A=x_{1}x_{2}\cdots x_{7}\in \C_{7}$, the vertex $x_{4}$ is called the {\it center} of $A$.
A directed path $\P=A_{1}A_{2}\cdots A_{l}~(l\geq 2)$ of $\D_{F}$ is {\it admissible} if $A_{1}\in \C(F)-(\C_{3}\cup \C^{(1)}_{5})$ and $A_{i}\in \C^{(2)}_{3}\cup \C^{(1)}_{5}$ for every $i~(2\leq i\leq l-1)$.
An admissible path $\P=A_{1}A_{2}\cdots A_{l}$ of $\D_{F}$ is {\it weakly admissible} if either
\begin{enumerate}[{\bf (W1)}]
\item
$A_{1}\in \C^{(2)}_{5}$ and $|V(A_{1})\cap T_{2}|=1$, or
\item
$A_{1}\in \C^{(1)}_{7}$ and $\sigma _{F}(A_{1}A_{2})$ is the center of $A_{1}$.
\end{enumerate} 
An admissible path $\P$ of $\D_{F}$ is {\it strongly admissible} if $\P$ is not weakly admissible.

A {\it path system} with respect to $F$ is a sequence $(\P_{1},\ldots ,\P_{m})~(m\geq 0)$ of admissible paths such that 
\begin{enumerate}[{\bf (P1)}]
\item
for each $i~(1\leq i\leq m)$, when we write $\P_{i}=A_{1}A_{2}\cdots A_{l}$, $\{A_{j}\mid 1\leq j\leq l-1\}\cap (\bigcup _{1\leq j\leq i-1}V(\P_{j}))=\emptyset $ and $A_{l}\in \C^{(1)}_{3}\cup (\bigcup _{1\leq j\leq i-1}V(\P_{j}))$, and
\item
for each $i~(1\leq i\leq m-1)$, $\P_{i}$ is weakly admissible.
\end{enumerate}
A path system $(\P_{1},\ldots ,\P_{m})$ with respect to $F$ is {\it complete} if $m\geq 1$ and $\P_{m}$ is strongly admissible.

By straightforward calculations, we get the following claim (and we omit its proof).

\begin{claim}%%%%%%%%%%%%%%%%%%%%%%%%%%%%%%%%%%%%%%%%%%%%%%%%%%%%%%%%%%%%%%%%%%%%%%%%%%%%%%%%%%%%%%%%%%%%%%%%%%%%%%%%%%%
\label{cl2.0}
Let $F$ be an $S$-central path-factor of $G$.
Then the following hold.
\begin{enumerate}[{\upshape(i)}]
\item
For $A\in \C^{(1)}_{3}(F)$, $|V(A)\cap T_{1}|+\frac{2}{3}|V(A)\cap T_{2}|=2=\frac{4}{3}|V(A)\cap S|+\frac{2}{3}$.
\item
For $A\in \C^{(2)}_{3}(F)$, $|V(A)\cap T_{1}|+\frac{2}{3}|V(A)\cap T_{2}|\geq \frac{4}{3}|V(A)\cap S|$.
\item
For $A\in \C^{(1)}_{5}(F)$, $|V(A)\cap T_{1}|+\frac{2}{3}|V(A)\cap T_{2}|>\frac{4}{3}|V(A)\cap S|$.
\item
For $A\in \C^{(2)}_{5}(F)$ with $|V(A)\cap T_{2}|=1$, $|V(A)\cap T_{1}|+\frac{2}{3}|V(A)\cap T_{2}|=\frac{4}{3}|V(A)\cap S|$.
\item
For $A\in \C^{(1)}_{7}(F)$, $|V(A)\cap T_{1}|+\frac{2}{3}|V(A)\cap T_{2}|=\frac{4}{3}|V(A)\cap S|$.
\qed
\end{enumerate}
\end{claim}
%%%%%%%%%%%%%%%%%%%%%%%%%%%%%%%%%%%%%%%%%%%%%%%%%%%%%%%%%%%%%%%%%%%%%%%%%%%%%%%%%%%%%%%%%%%%%%%%%%%%%%%%%%%%%%%%%%%%%%%%

The following claim plays a key role in the proof of the theorem.

\begin{claim}%%%%%%%%%%%%%%%%%%%%%%%%%%%%%%%%%%%%%%%%%%%%%%%%%%%%%%%%%%%%%%%%%%%%%%%%%%%%%%%%%%%%%%%%%%%%%%%%%%%%%%%%%%%
\label{cl2.1}
Let $F$ be an $S$-central path-factor of $G$ with $\C^{(1)}_{3}(F)\not=\emptyset $, and let $(\P_{1},\ldots ,\P_{m})$ be a path system with respect to $F$ ($m\geq 0$).
Then the system can be extend to a complete path system $(\P_{1},\ldots ,\P_{m},\P_{m+1},\ldots ,\P_{m'})$ with respect to $F$.
\end{claim}
%%%%%%%%%%%%%%%%%%%%%%%%%%%%%%%%%%%%%%%%%%%%%%%%%%%%%%%%%%%%%%%%%%%%%%%%%%%%%%%%%%%%%%%%%%%%%%%%%%%%%%%%%%%%%%%%%%%%%%%%
\proof
We take a maximal path system $(\P_{1},\ldots ,\P_{m},\P_{m+1},\ldots ,\P_{m'})$ with respect to $F$.
We show that $(\P_{1},\ldots ,\P_{m'})$ is a complete path system.
Suppose that $(\P_{1},\ldots ,\P_{m'})$ is not a complete path system.
Then $\P_{i}$ is weakly admissible for each $i$ with $1\leq i\leq m'$ (this includes the case where $m'=0$).

Set $\A_{1}=\bigcup _{1\leq i\leq m'}V(\P_{i})$ (note that $\A_{1}=\emptyset $ if and only if $m'=0$).
Let $X=(\bigcup _{A\in \A_{1}}V(A))\cap S$ and $Y_{h}=(\bigcup _{A\in \A_{1}}V(A))\cap T_{h}~(h\in \{1,2\})$.
Then by the definition of a weakly admissible path (and the definition of a path system), $\A_{1}\subseteq \C_{3}\cup \C_{5}\cup \C^{(1)}_{7}$, and if $A\in \A_{1}\cap \C^{(2)}_{5}$, then $|V(A)\cap T_{1}|=1$.
Furthermore, by condition (P1) in the definition of a path system, $\A_{1}\not=\emptyset $ if and only if $\A_{1}\cap \C^{(1)}_{3}\not=\emptyset $.
Hence by Claim~\ref{cl2.0},
\begin{align}
|Y_{1}|+\frac{2}{3}|Y_{2}|\geq \frac{4}{3}|X|\label{cl2.1-6}
\end{align}
and
\begin{align}
|Y_{1}|+\frac{2}{3}|Y_{2}|\geq \frac{4}{3}|X|+\frac{2}{3} \mbox{ if } \A_{1}\not=\emptyset .\label{cl2.1-7}
\end{align}

Let $\A_{2}=\C^{(1)}_{3}-\A_{1}$, $X^{*}=(\bigcup _{A\in \A_{2}}V(A))\cap S$ and $Y^{*}_{h}=(\bigcup _{A\in \A_{2}}V(A))\cap T_{h}~(h\in \{1,2\})$.
By Claim~\ref{cl2.0}(i),
\begin{align}
|Y^{*}_{1}|+\frac{2}{3}|Y^{*}_{2}|\geq \frac{4}{3}|X^{*}|\label{cl2.1-8}
\end{align}
and
\begin{align}
|Y^{*}_{1}|+\frac{2}{3}|Y^{*}_{2}|=\frac{4}{3}|X^{*}|+\frac{2}{3} \mbox{ if } \A_{2}\not=\emptyset .\label{cl2.1-9}
\end{align}

Let $(B_{1},\ldots ,B_{l})~(l\geq 0)$ be a sequence such that for each $i~(1\leq i\leq l)$, $B_{i}\in (\C^{(2)}_{3}\cup \C^{(1)}_{5})-(\A_{1}\cup \A_{2}\cup \{B_{j}\mid 1\leq j\leq i-1\})$ and there exists an edge of $\D_{F}$ from $B_{i}$ to an element in $\A_{1}\cup \A_{2}\cup \{B_{j}\mid 1\leq j\leq i-1\}$.
We choose $(B_{1},\ldots ,B_{l})$ so that $l$ is as large as possible.
Let $\A_{3}=\{B_{i}\mid 1\leq i\leq l\}$, $X^{**}=(\bigcup _{A\in \A_{3}}V(A))\cap S$ and $Y^{**}_{h}=(\bigcup _{A\in \A_{3}}V(A))\cap T_{h}~(h\in \{1,2\})$.
By Claim~\ref{cl2.0}(ii)(iii),
\begin{align}
|Y^{**}_{1}|+\frac{2}{3}|Y^{**}_{2}|\geq \frac{4}{3}|X^{**}|.\label{cl2.1-10}
\end{align}

Let $X^{0}=X\cup X^{*}\cup X^{**}$ and $Y^{0}_{h}=Y_{h}\cup Y^{*}_{h}\cup Y^{**}_{h}~(h\in \{1,2\})$.
If $m'\geq 1$, then $\A_{1}\not=\emptyset $; if $m'=0$ (i.e., $\A_{1}=\emptyset $), then $\A_{2}\not=\emptyset $ because $\C^{(1)}_{3}\not=\emptyset $.
Thus by (\ref{cl2.1-7}) and (\ref{cl2.1-9}), either $|Y_{1}|+\frac{2}{3}|Y_{2}|\geq \frac{4}{3}|X|+\frac{2}{3}$ or $|Y^{*}_{1}|+\frac{2}{3}|Y^{*}_{2}|=\frac{4}{3}|X^{*}|+\frac{2}{3}$.
This together with (\ref{cl2.1-6}), (\ref{cl2.1-8}) and (\ref{cl2.1-10}) leads to
\begin{align}
|Y^{0}_{1}|+\frac{2}{3}|Y^{0}_{2}|\geq \frac{4}{3}|X^{0}|+\frac{2}{3}.\label{cl2.1-11}
\end{align}
Since $|T_{1}|+\frac{2}{3}|T_{2}|\leq \frac{4}{3}|S|+\frac{1}{3}$, this implies $X^{0}\not=S$ and hence $\C(F)-(\A_{1}\cup \A_{2}\cup \A_{3})\not=\emptyset $.

Let $\tilde{\A}=\C(F)-(\A_{1}\cup \A_{2}\cup \A_{3})$, $\tilde{X}=(\bigcup _{A\in \tilde{\A}}V(A))\cap S$ and $\tilde{Y}_{h}=(\bigcup _{A\in \tilde{\A}}V(A))\cap T_{h}~(h\in \{1,2\})$.
Note that $S$ is the disjoint union of $X^{0}$ and $\tilde{X}$ and, for $h\in \{1,2\}$, $T_{h}$ is the disjoint union of $Y^{0}_{h}$ and $\tilde{Y}_{h}$.
If $|\tilde{Y}_{1}|+\frac{2}{3}|\tilde{Y}_{2}|\geq \frac{4}{3}|\tilde{X}|$, then by (\ref{cl2.1-11}), $|T_{1}|+\frac{2}{3}|T_{2}|=(|Y^{0}_{1}|+|\tilde{Y}_{1}|)+\frac{2}{3}(|Y^{0}_{2}|+|\tilde{Y}_{2}|)\geq \frac{4}{3}|X^{0}|+\frac{2}{3}+\frac{4}{3}|\tilde{X}|=\frac{4}{3}|S|+\frac{2}{3}$, which is a contradiction.
Thus $|\tilde{Y}_{1}|+\frac{2}{3}|\tilde{Y}_{2}|<\frac{4}{3}|\tilde{X}|$.
On the other hand, since $\A_{1}\cup \A_{2}\not=\emptyset $, we have $Y^{0}_{1}\cup Y^{0}_{2}\not=\emptyset $, and hence $\tilde{Y}_{1}\cup \tilde{Y}_{2}\not=T$.
Consequently $N_{G}(\tilde{X})\not\subseteq \tilde{Y}_{1}\cup \tilde{Y}_{2}$ by the assumption of the theorem, which implies that there exists a vertex $x\in \tilde{X}$ with $N_{G}(x)\cap (Y^{0}_{1}\cup Y^{0}_{2})\not=\emptyset $.
Let $\tilde{A}\in \tilde{\A}$ be the path containing $x$.
By the definition of $\A_{2}$ and $\tilde{\A}$, $\tilde{A}\not\in \C^{(1)}_{3}$.
By the maximality of $(B_{1},\ldots ,B_{l})$, $\tilde{A}\not\in \C^{(2)}_{3}\cup \C^{(1)}_{5}$.
Thus $\tilde{A}\in \C(F)-(\C_{3}\cup \C^{(1)}_{5})$.
By the definition of $(B_{1},\ldots ,B_{l})$ and $x$, there exists a directed path $\P'=\tilde{A}_{1}\cdots \tilde{A}_{p}$ of $\D_{F}$ such that $\tilde{A}_{1}=\tilde{A}$, $\tilde{A}_{i}\in \A_{3}~(2\leq i\leq p-1)$ and $\tilde{A}_{p}\in \A_{1}\cup \A_{2}$.
Then $\P'$ is an admissible path of $\D_{F}$.
Now the sequence $(\P_{1},\ldots ,\P_{m'},\P')$ is a path system with respect to $F$, which contradicts the maximality of $(\P_{1},\ldots ,\P_{m'})$.
This contradiction completes the proof of the claim.
\qed

We turn to the proof of Theorem~\ref{thm2}.
By way of contradiction, suppose that $\C^{(1)}_{3}(F)\not=\emptyset $ for every $S$-central path-factor $F$ of $G$.
By Lemma~\ref{lem2.1}, $G$ has an $S$-central path-factor $F_{0}$.
Note that an empty sequence is a path system with respect to $F_{0}$.
Hence by Claim~\ref{cl2.1}, there exists a complete path system $(\P_{1},\ldots ,\P_{m})$ with respect to $F_{0}$.
Choose $F_{0}$ and $(\P_{1},\ldots ,\P_{m})$ so that
\begin{enumerate}[{\bf (F1)}]
\item
$|\C^{(1)}_{3}(F_{0})|$ is as small as possible, and
\item
subject to (F1), $(|V(\P_{1})|,\ldots ,|V(\P_{m})|)$ is lexicographically as small as possible.
\end{enumerate}

For each $i~(1\leq i\leq m)$, write $\P_{i}=A^{(i)}_{1}\cdots A^{(i)}_{l_{i}}$.
Then $\bigcup _{1\leq i\leq m}\P_{i}$ contains a directed path $B_{1}B_{2}\cdots B_{p}$ of $\D_{F_{0}}$ with $B_{1}=A^{(m)}_{1}$ and $B_{p}\in \C^{(1)}_{3}(F_{0})$.
For each $i~(1\leq i\leq p)$, write $B_{i}=v_{i,1}v_{i,2}\cdots v_{i,q_{i}}$.
For $i~(1\leq i\leq p-1)$, let $s_{i}$ be the integer with $v_{i,s_{i}}=\sigma _{F}(B_{i}B_{i+1})$, and for $i~(2\leq i\leq p)$, let $t_{i}$ be the integer with $v_{i,t_{i}}=\tau _{F_{0}}(B_{i-1}B_{i})$.
As in the proof of Claim~\ref{cl2.1.1}, by renumbering the vertices of some of the $B_{i}$ backward if necessary, we may assume that
\begin{enumerate}[{\bf (B1)}]
\item
$s_{1}\geq \frac{q_{1}+1}{2}$ if $q_{1}$ is odd,
\item
$\{v_{1,1},v_{1,3}\}\cap T_{2}\not=\emptyset $ if $B_{1}\in \C^{(2)}_{7}(F_{0})$ and $s_{1}=4$,
\item
$s_{1}$ is odd if $q_{1}$ is even,
\item
$t_{i}<s_{i}$ for each $i~(2\leq i\leq p-1)$, and
\item
$t_{p}=q_{p}~(=3)$.
\end{enumerate}
Note that (B3) means that when $q_{1}$ is even, the vertices of $B_{1}$ are numbered so that $v_{1,q_{1}}\in T$.
Thus $v_{i,q_{i}}\in T$ for each $i~(1\leq i\leq p)$.
We can divide the type of $B_{1}$ into three possibilities as follows:

\begin{claim}%%%%%%%%%%%%%%%%%%%%%%%%%%%%%%%%%%%%%%%%%%%%%%%%%%%%%%%%%%%%%%%%%%%%%%%%%%%%%%%%%%%%%%%%%%%%%%%%%%%%%%%%%%%
\label{cl-type1}
One of the following holds:
\begin{enumerate}[{\upshape(1)}]
\item
$|V(B_{1})|$ is even and $s_{1}$ is odd;
\item
$B_{1}\in \C^{(2)}_{5}(F_{0})\cup \C^{(2)}_{7}(F_{0})$, $s_{1}=4$ and $\{v_{1,1},v_{1,3}\}\cap T_{2}\not=\emptyset $; or
\item
$|V(B_{1})|\geq 7$ and $s_{1}\geq 6$.
\end{enumerate}
\end{claim}
%%%%%%%%%%%%%%%%%%%%%%%%%%%%%%%%%%%%%%%%%%%%%%%%%%%%%%%%%%%%%%%%%%%%%%%%%%%%%%%%%%%%%%%%%%%%%%%%%%%%%%%%%%%%%%%%%%%%%%%%
\proof
If $|V(B_{1})|$ is even, then (1) holds by (B3).
Thus we may assume $|V(B_{1})|$ is odd.
Then by the definition of a strongly admissible path, $B_{1}\in \C^{(2)}_{5}(F_{0})$ and $|V(B_{1})\cap T_{2}|\geq 2$, or $B_{1}\in \C^{(1)}_{7}(F_{0})$ and $s_{1}\not=4$, or $B_{1}\in \C^{(2)}_{7}(F_{0})$, or $|V(B_{1})|\geq 9$.
If $B_{1}\in \C^{(2)}_{5}(F_{0})$ and $|V(B_{1})\cap T_{2}|\geq 2$, then (2) holds by (B1).
If $B_{1}\in \C^{(1)}_{7}(F_{0})\cup \C^{(2)}_{7}(F_{0})$ and $s_{1}\not=4$, then (3) holds by (B1).
If $B_{1}\in \C^{(2)}_{7}(F_{0})$ and $s_{1}=4$, then (2) holds by (B2).
If $|V(B_{1})|\geq 9$, then (3) holds by (B1).
\qed

As for $B_{i}$ with $2\leq i\leq p-1$, the following claim follows immediately from the definition of a weakly admissible path.

\begin{claim}%%%%%%%%%%%%%%%%%%%%%%%%%%%%%%%%%%%%%%%%%%%%%%%%%%%%%%%%%%%%%%%%%%%%%%%%%%%%%%%%%%%%%%%%%%%%%%%%%%%%%%%%%%%
\label{cl-type2}
Let $2\leq i\leq p-1$.
Then one of the following holds:
\begin{enumerate}[{\upshape(1)}]
\item
$B_{i}\in \C^{(2)}_{3}(F_{0})$ and $s_{i}=2$;
\item
$B_{i}\in \C_{5}(F_{0})$ and $s_{i}=2$ or $4$; or
\item
$B_{i}\in \C^{(1)}_{7}(F_{0})$ and $s_{i}=4$.
\qed
\end{enumerate}
\end{claim}
%%%%%%%%%%%%%%%%%%%%%%%%%%%%%%%%%%%%%%%%%%%%%%%%%%%%%%%%%%%%%%%%%%%%%%%%%%%%%%%%%%%%%%%%%%%%%%%%%%%%%%%%%%%%%%%%%%%%%%%%

Let $i_{0}$ be the minimum integer $i~(\geq 2)$ satisfying one of the following two conditions:
\begin{enumerate}
\item[{\bf (I1)}]
$i=p$; or
\item[{\bf (I2)}]
$2\leq i\leq p-1$ and $t_{i}=1$.
\end{enumerate}
Set $B'_{1}=v_{1,1}v_{1,2}\cdots v_{1,s_{1}-1}$ and, for each $i~(2\leq i\leq i_{0})$, set
$$
B'_{i}=v_{i-1,q_{i-1}}v_{i-1,q_{i-1}-1}\cdots v_{i-1,s_{i-1}}v_{i,t_{i}}v_{i,t_{i}-1}\cdots v_{i,1}
$$
(see Figure~\ref{thm2.1fig2}).
Let $2\leq i\leq i_{0}-1$.
By the definition of $i_{0}$, $t_{i}\geq 3$.
On the other hand, $s_{i}\leq 4$ by Claim~\ref{cl-type2}.
Hence $t_{i}=s_{i}-1$.
Since $i~(2\leq i\leq i_{0}-1)$ is arbitrary, it follows that
\begin{align}
\label{neweq2.1}
B'_{1},\ldots ,B'_{i_{0}}\mbox{ are vertex-disjoint paths of }G
\end{align}
and
\begin{align}
\label{neweq2.2}
\bigcup _{1\leq i\leq i_{0}}V(B'_{i})=\bigcup _{1\leq i\leq i_{0}}V(B_{i})-\{v_{i_{0},j}\mid t_{i_{0}}+1\leq j\leq q_{i_{0}}\}.
\end{align}
Furthermore, 
\begin{align}
\label{neweq2.3}
|V(B'_{i})\cap T|\geq |V(B'_{i})\cap S|\mbox{ for each }i~(2\leq i\leq i_{0})
\end{align}
because $v_{i-1,q_{i-1}}\in T$.
If $B'_{1}\not=\emptyset $, then $v_{1,s_{1}-1}\in T$, and hence
\begin{align}
\label{neweq2.4}
|V(B'_{1})\cap T|\geq |V(B'_{1})\cap S|
\end{align}
(if $B'_{1}=\emptyset $, then (\ref{neweq2.4}) trivially holds).
Also
\begin{align}
\label{neweq2.5}
|V(B'_{i})\cap V(B_{i-1})|\mbox{ is even and }|V(B'_{i})\cap V(B_{i-1})|\geq 2\mbox{ for each }i~(2\leq i\leq i_{0})
\end{align}
because $v_{i-1,s_{i-1}}\in S$ and $v_{i-1,q_{i-1}}\in T$.
It follows from (\ref{neweq2.5}) that
\begin{align}
\label{neweq2.6}
|V(B'_{i})|\geq 5\mbox{ for each }i~(2\leq i\leq i_{0}-1)
\end{align}
because $|V(B'_{i})\cap V(B_{i})|=t_{i}\geq 3$.
Since $|V(B'_{1})|=s_{1}-1$, we see from Claim~\ref{cl-type1} that
\begin{align}
\label{neweq2.7}
|V(B'_{1})|\mbox{ is even or }|V(B'_{1})|\geq 3,
\end{align}
and
\begin{align}
\label{neweq2.8}
V(B'_{1})\cap T_{2}\not=\emptyset \mbox{ if }|V(B'_{1})|=3.
\end{align}
Combining (\ref{neweq2.3}) through (\ref{neweq2.8}), we get the following claim.

\begin{figure}
\begin{center}
%WinTpicVersion4.26
\unitlength 0.1in
%
\caption{Construction of $B'_{i}$}
\label{thm2.1fig2}
\end{center}
\end{figure}

\begin{claim}%%%%%%%%%%%%%%%%%%%%%%%%%%%%%%%%%%%%%%%%%%%%%%%%%%%%%%%%%%%%%%%%%%%%%%%%%%%%%%%%%%%%%%%%%%%%%%%%%%%%%%%%%%%
\label{cl2.2}
\begin{enumerate}[{\upshape(i)}]
\item
For each $i$ with $1\leq i\leq i_{0}$, we have $|V(B'_{i})\cap T|\geq |V(B'_{i})\cap S|$.
\item
For each $i$ with $1\leq i\leq i_{0}-1$,
\begin{enumerate}[{\upshape(a)}]
\item
$|V(B'_{i})|$ is even or $|V(B'_{i})|\geq 3$, and
\item
$V(B'_{i})\cap T_{2}\not=\emptyset $ if $B'_{i}\simeq P_{3}$.
\qed
\end{enumerate}
\end{enumerate}
\end{claim}
%%%%%%%%%%%%%%%%%%%%%%%%%%%%%%%%%%%%%%%%%%%%%%%%%%%%%%%%%%%%%%%%%%%%%%%%%%%%%%%%%%%%%%%%%%%%%%%%%%%%%%%%%%%%%%%%%%%%%%%%

Suppose that $i_{0}=p$.
Then
\begin{align}
\label{neweq2.9}
|V(B'_{p})|\geq 5
\end{align}
by (\ref{neweq2.5}) and (B5).
Let $F_{1}=(F_{0}-(\bigcup _{1\leq i\leq p}V(B_{i})))\cup (\bigcup _{1\leq i\leq p}B'_{i})$.
Then by Claim~\ref{cl2.2}, (\ref{neweq2.9}), (\ref{neweq2.1}), (\ref{neweq2.2}) and (B5), $F_{1}$ is an $S$-central path-factor of $G$, and $B'_{i}\not\in \C^{(1)}_{3}(F_{1})$ for each $i~(1\leq i\leq p)$.
Since $B_{1}\in \C^{(1)}_{3}(F_{0})$, we have $|\C^{(1)}_{3}(F_{1})|<|\C^{(1)}_{3}(F_{0})|$, which contradicts the minimality of $|\C^{(1)}_{3}(F)|$.
Thus $2\leq i_{0}\leq p-1$.
Then by the definition of $i_{0}$, $t_{i_{0}}=1$.
Hence $B''_{i_{0}}=B_{i_{0}}\cup B'_{i_{0}}$ is a path of $G$ with $|V(B''_{i_{0}})\cap T|\geq |V(B''_{i_{0}})\cap S|$ (see Figure~\ref{thm2.1fig3}).
Set $F_{2}=(F_{0}-(\bigcup _{1\leq i\leq i_{0}}V(B_{i})))\cup (\bigcup _{1\leq i\leq i_{0}-1}B'_{i})\cup B''_{i_{0}}$.
Then by Claim~\ref{cl2.2}, (\ref{neweq2.1}) and (\ref{neweq2.2}), $F_{2}$ is an $S$-central path-factor of $G$, and $B'_{i}\not\in \C^{(1)}_{3}(F_{1})$ for each $i~(1\leq i\leq i_{0}-1)$.
Furthermore,
\begin{align}
|V(B''_{i_{0}})|=|V(B_{i_{0}})|+|V(B'_{i_{0}})\cap V(B_{i_{0}-1})|.\label{cl2.1-12}
\end{align}
Since $|V(B'_{i_{0}})\cap V(B_{i_{0}-1})|\geq 2$ by (\ref{neweq2.5}), this implies $|V(B''_{i_{0}})|\geq 5$, and hence we also have $B''_{i_{0}}\not\in \C^{(1)}_{3}(F_{1})$.
Thus $|\C^{(1)}_{3}(F_{2})|=|\C^{(1)}_{3}(F_{0})|$.

\begin{figure}
\begin{center}
%WinTpicVersion4.26
\unitlength 0.1in
%
\caption{Construction of $B''_{i_{0}}$}
\label{thm2.1fig3}
\end{center}
\end{figure}

Set $k_{0}=\min\{k\mid B_{i_{0}}\in V(\P_{k})\}$, and write $B_{i_{0}}=A^{(k_{0})}_{j_{0}}$.
If $B_{p}\in V(\P_{k_{0}})$, then the fact that $B_{i_{0}}\not=B_{p}$ implies that $j_{0}\leq l_{k_{0}}-1$; if $B_{p}\not\in V(\P_{k_{0}})$, then the minimality of $k_{0}$ implies that $j_{0}\leq l_{k_{0}}-1$.
In either case, we have $j_{0}\leq l_{k_{0}}-1$.

\medskip
\noindent
\textbf{Case 1:} $j_{0}=1$.

Since $B_{i_{0}}=A^{(k_{0})}_{1}$ and $i_{0}\geq 2$, $B_{1}\in \bigcup _{k_{0}+1\leq i\leq m}V(\P_{i})$.
In particular, $k_{0}\leq m-1$ and $\P_{k_{0}}$ is weakly admissible.
Hence $B_{i_{0}}\in \C^{(2)}_{5}(F_{0})\cup \C^{(1)}_{7}(F_{0})$.
This together with (\ref{cl2.1-12}) and (\ref{neweq2.5}) implies that $B''_{i_{0}}\in \C^{(2)}_{7}(F_{2})$ or $|V(B''_{i_{0}})|\geq 9$.
Thus the directed path $\P'_{k_{0}}=B''_{i_{0}}A^{(k_{0})}_{2}\cdots A^{(k_{0})}_{l_{k_{0}}}$ of $\D_{F_{2}}$ is strongly admissible.
Consequently $(\P_{1},\ldots ,\P_{k_{0}-1},\P'_{k_{0}})$ is a complete path system with respect to $F_{2}$.
Since $k_{0}\leq m-1$ and $|V(\P_{k_{0}})|=|V(\P'_{k_{0}})|$, we see that $(|V(\P_{1})|,\ldots ,|V(\P_{k_{0}-1})|,|V(\P'_{k_{0}})|)$ is lexicographically less than $(|V(\P_{1})|,\ldots ,|V(\P_{k_{0}-1})|,|V(\P_{k_{0}})|,\ldots ,|V(\P_{m})|)$, which contradicts the minimality of $(|V(\P_{1})|,\ldots ,|V(\P_{m})|)$.

\medskip
\noindent
\textbf{Case 2:} $2\leq j_{0}\leq l_{k_{0}}-1$.

Since $B_{i_{0}}=A^{(k_{0})}_{j_{0}}$, $B_{i_{0}}\in \C^{(2)}_{3}(F_{0})\cup \C^{(1)}_{5}(F_{0})$.
This together with (\ref{cl2.1-12}) and (\ref{neweq2.5}) implies that $B''_{i_{0}}\in \C^{(2)}_{5}(F_{2})$ or $|V(B''_{i_{0}})|\geq 7$.
Thus the directed path $\P'_{k_{0}}=B''_{i_{0}}A^{(k_{0})}_{j_{0}+1}A^{(k_{0})}_{j_{0}+2}\cdots A^{(k_{0})}_{l_{k_{0}}}$ of $\D_{F_{2}}$ is admissible.
Consequently $(\P_{1},\ldots ,\P_{k_{0}-1},\P'_{k_{0}})$ is a path system with respect to $F_{2}$.
By Claim~\ref{cl2.1}, the system can be extend to a complete path system $(\P_{1},\ldots ,\P_{k_{0}-1},\P'_{k_{0}},\Q_{1},\ldots ,\Q_{\alpha })$ with respect to $F_{2}$ (it is possible that $\alpha =0$).
Since $j_{0}\geq 2$, $|V(\P'_{k_{0}})|=l_{k_{0}}-j_{0}+1<l_{k_{0}}=|V(\P_{k_{0}})|$, and hence $(|V(\P_{1})|,\ldots ,|V(\P_{k_{0}-1})|,|V(\P'_{k_{0}})|,|V(\Q_{1})|,\ldots ,|V(\Q_{\alpha })|)$ is lexicographically less than $(|V(\P_{1})|,\ldots ,|V(\P_{k_{0}-1})|,|V(\P_{k_{0}})|,\ldots ,|V(\P_{m})|)$, which contradicts the minimality of $(|V(\P_{1})|,\ldots ,|V(\P_{m})|)$.

This completes the proof of Theorem~\ref{thm2}.
\qed

%%%%%%%%%%%%%%%%%%%%%%%%%%%%%%%%%%%%%%%%%%%%%%%%%%%%%%%%%%%%%%%%%%%%%%%%%%%%%%%%%%%%%%%%%%%%%%%%%%%%%%%%%%%%%%%%%%%%%%%%
%%%%%%%%%%%%%%%%%%%%%%%%%%%%%%%%%%%%%%%%%%%%%%%%%%%%%%%%%%%%%%%%%%%%%%%%%%%%%%%%%%%%%%%%%%%%%%%%%%%%%%%%%%%%%%%%%%%%%%%%
\section{Proof of Theorem~\ref{thm1}}\label{sec3}
%%%%%%%%%%%%%%%%%%%%%%%%%%%%%%%%%%%%%%%%%%%%%%%%%%%%%%%%%%%%%%%%%%%%%%%%%%%%%%%%%%%%%%%%%%%%%%%%%%%%%%%%%%%%%%%%%%%%%%%%
%%%%%%%%%%%%%%%%%%%%%%%%%%%%%%%%%%%%%%%%%%%%%%%%%%%%%%%%%%%%%%%%%%%%%%%%%%%%%%%%%%%%%%%%%%%%%%%%%%%%%%%%%%%%%%%%%%%%%%%%

Let $G$ be as in Theorem~\ref{thm1}.
By assumption, we have $c_{1}(G)+\frac{2}{3}c_{3}(G)\leq \frac{4}{3}|\emptyset |+\frac{1}{3}=\frac{1}{3}$.
Hence $c_{1}(G)=c_{3}(G)=0$.

We now proceed by induction on $|V(G)|+|E(G)|$.
We may assume $V(G)\not=\emptyset $.
Note that if $E(G)=\emptyset $, then $c_{1}(G)=|V(G)|\geq 1$, which is a contradiction.
This means that the theorem holds for graphs $G$ with $E(G)=\emptyset $ in the sense that the assumption is not satisfied.
We henceforth assume that $E(G)\not=\emptyset $ and the theorem holds for graphs $G'$ with $|V(G')|+|E(G')|<|V(G)|+|E(G)|$.

Let $\S =\{X\subseteq V(G)\mid c_{1}(G-X)+c_{3}(G-X)\geq 1\}$.
Since $c_{1}(G-N_{G}(x))\geq 1$ for $x\in V(G)$, $\S \not=\emptyset $.
Set 
$$
\beta =\min_{X\in \S}\left\{\frac{4}{3}|X|+\frac{1}{3}-c_{1}(G-X)-\frac{2}{3}c_{3}(G-X)\right\}.
$$
%By the assumption of the theorem, $\beta \geq 0$.

\begin{claim}%%%%%%%%%%%%%%%%%%%%%%%%%%%%%%%%%%%%%%%%%%%%%%%%%%%%%%%%%%%%%%%%%%%%%%%%%%%%%%%%%%%%%%%%%%%%%%%%%%%%%%%%%%%
\label{cl3.0}
If $\beta \geq 2$, then $G$ has a $\{P_{2},P_{5}\}$-factor.
\end{claim}
%%%%%%%%%%%%%%%%%%%%%%%%%%%%%%%%%%%%%%%%%%%%%%%%%%%%%%%%%%%%%%%%%%%%%%%%%%%%%%%%%%%%%%%%%%%%%%%%%%%%%%%%%%%%%%%%%%%%%%%%
\proof
Let $e\in E(G)$, and suppose that $\C_{1}(G-e)\cup \C_{3}(G-e)\not=\emptyset $.
Take $C\in \C_{1}(G-e)\cup \C_{3}(G-e)$.
Since $c_{1}(G)=c_{3}(G)=0$, $e$ joins a vertex in $V(C)$ and a vertex $y$ in $V(G)-V(C)$.
This implies $C\in \C_{1}(G-y)\cup \C_{3}(G-y)$, and hence $\frac{4}{3}|\{y\}|+\frac{1}{3}-(c_{1}(G-y)+\frac{2}{3}c_{3}(G-y))\leq \frac{4}{3}+\frac{1}{3}-\frac{2}{3}=1$, which contradicts the assumption that $\beta \geq 2$.
Thus $c_{1}(G-e)=c_{3}(G-e)=0$ for all $e\in E(G)$.
From the fact that $c_{1}(G-e)=\emptyset $ for all $e\in E(G)$, it follows that $d_{G}(x)\geq 2$ for all $x\in V(G)$.
Assume for the moment that $d_{G}(x)=2$ for all $x\in V(G)$.
Then each component of $G$ is a cycle.
Since $c_{3}(G)=0$, this implies that $G$ has a path-factor $F$ with $\C_{3}(F)=\emptyset $.
Hence by Fact~\ref{fact1}, $G$ has a $\{P_{2},P_{5}\}$-factor.
Thus we may assume that there exists $x_{0}\in V(G)$ such that $d_{G}(x_{0})\geq 3$.

Fix an edge $e^{*}=x_{0}y_{0}\in E(G)$ incident with $x_{0}$, and let $G'=G-e^{*}$.
By an assertion in the first paragraph of the proof the claim, $c_{1}(G')=c_{3}(G')=0$.
Let $X\subseteq V(G')$.
We show that $\frac{4}{3}|X|+\frac{1}{3}-c_{1}(G'-X)-\frac{2}{3}c_{3}(G'-X)\geq 0$.
We have $\frac{4}{3}|\emptyset |+\frac{1}{3}-c_{1}(G')-c_{3}(G')=\frac{1}{3}>0$.
Thus we may assume $X\not=\emptyset $.
Note that $|(\C_{1}(G'-X)-\C_{1}(G-X)|+|(\C_{3}(G'-X)-\C_{3}(G-X)|\leq 2$, and hence
\begin{align}
c_{1}(G'-X)+\frac{2}{3}c_{3}(G'-X)\leq c_{1}(G-X)+\frac{2}{3}c_{3}(G-X)+2.\label{cl2.0-00}
\end{align}
Furthermore, if equality holds in (\ref{cl2.0-00}), then $x_{0},y_{0}\not\in X$ and $\{x_{0}\},\{y_{0}\}\in \C_{1}(G'-X)$.
If $c_{1}(G-X)+c_{3}(G-X)\geq 1$, then by the definition of $\beta $, $\frac{4}{3}|X|+\frac{1}{3}-c_{1}(G-X)-\frac{2}{3}c_{3}(G-X)\geq \beta \geq 2$ which, together with (\ref{cl2.0-00}), leads to $\frac{4}{3}|X|+\frac{1}{3}-(c_{1}(G'-X)+\frac{2}{3}c_{3}(G'-X))\geq \frac{4}{3}|X|+\frac{1}{3}-(c_{1}(G-X)+\frac{2}{3}c_{3}(G-X)+2)\geq \beta -2\geq 0$.
Thus we may assume that $c_{1}(G-X)+c_{3}(G-X)=0$.
By (\ref{cl2.0-00}), $c_{1}(G'-X)+c_{3}(G'-X)\leq c_{1}(G-X)+c_{3}(G-X)+2=2$.
By way of contradiction, suppose that $\frac{4}{3}|X|+\frac{1}{3}-(c_{1}(G'-X)+\frac{2}{3}c_{3}(G'-X))<0$.
Then $\frac{4}{3}|X|+\frac{1}{3}-2<0$.
Since $X\not=\emptyset $, this forces $|X|=1$ and $c_{1}(G'-X)+\frac{2}{3}c_{3}(G'-X)=2$.
Hence equality in (\ref{cl2.0-00}), which implies $\{x_{0}\}\in \C_{1}(G'-X)$.
Consequently $d_{G}(x_{0})\leq |X\cup \{y_{0}\}|=2$, which contradicts the fact that $d_{G}(x_{0})\geq 3$.
Thus we have $\frac{4}{3}|X|+\frac{1}{3}-c_{1}(G'-X)-\frac{2}{3}c_{3}(G'-X)\geq 0$ for all $X\subseteq V(G')$.
By the induction assumption, $G'$ has a $\{P_{2},P_{5}\}$-factor.
Therefore $G$ also has a $\{P_{2},P_{5}\}$-factor.
\qed

By Claim~\ref{cl3.0}, we may assume that $\beta \leq \frac{5}{3}$.

Let $S\in \S $ be a maximum set with $\frac{4}{3}|S|-c_{1}(G-S)-\frac{2}{3}c_{3}(G-S)+\frac{1}{3}=\beta $.
%Choose $S$ so that $|S|$ is as small as possible.

\begin{claim}%%%%%%%%%%%%%%%%%%%%%%%%%%%%%%%%%%%%%%%%%%%%%%%%%%%%%%%%%%%%%%%%%%%%%%%%%%%%%%%%%%%%%%%%%%%%%%%%%%%%%%%%%%%
\label{cl3.1}
Let $C$ be a component of $G-S$.
\begin{enumerate}[{\upshape(i)}]
\item If $|V(C)|\not\in \{1,3\}$, then $C$ has a $\{P_{2},P_{5}\}$-factor.
\item If $|V(C)|=3$, then $C$ is complete.
\end{enumerate}
\end{claim}
%%%%%%%%%%%%%%%%%%%%%%%%%%%%%%%%%%%%%%%%%%%%%%%%%%%%%%%%%%%%%%%%%%%%%%%%%%%%%%%%%%%%%%%%%%%%%%%%%%%%%%%%%%%%%%%%%%%%%%%%
\proof
\begin{enumerate}[{\upshape(i)}]
\item
Suppose that $C$ has no $\{P_{2},P_{5}\}$-factor.
Then by the induction assumption, there exists a set $S'\subseteq V(C)$ with $\frac{4}{3}|S'|+\frac{1}{3}-c_{1}(C-S')-\frac{2}{3}c_{3}(C-S')<0$.
Set $S_{0}=S\cup S'$.
Since $\C_{1}(G-S_{0})=\C_{1}(G-S)\cup \C_{1}(C-S')$, $\C_{3}(G-S_{0})=\C_{3}(G-S)\cup \C_{3}(C-S')$ and $\C_{1}(C-S')\cup \C_{3}(C-S')\not=\emptyset $, we have $S_{0}\in \S $.
We also get $\frac{4}{3}|S_{0}|+\frac{1}{3}-c_{1}(G-S_{0})-\frac{2}{3}c_{3}(G-S_{0})=(\frac{4}{3}|S|+\frac{1}{3}-c_{1}(G-S)-\frac{2}{3}c_{3}(G-S))+(\frac{4}{3}|S'|-c_{1}(C-S')-\frac{2}{3}c_{3}(C-S'))<\beta  $.
This contradicts the definition of $\beta $.
\item
Suppose that $|V(C)|=3$ and $C$ is not complete (i.e., $C$ is a path of order three).
Let $x\in C$ be the vertex with $d_{C}(x)=2$.
Then $c_{1}(C-x)=2$ and $c_{3}(C-x)=0$.
Set $S_{1}=S\cup \{x\}$.
%Then $\frac{4}{3}|\{x\}|-c_{1}(C-x)-\frac{2}{3}c_{3}(C-x)=\frac{4}{3}-2=-\frac{2}{3}$.
Since $\C_{1}(G-S_{1})=\C_{1}(G-S)\cup \C_{1}(C-x)$, $\C_{3}(G-S_{1})=\C_{3}(G-S)-\{C\}$ and $\C_{1}(C-x)\not=\emptyset $, we have $S_{1}\in \S $.
We also get $\frac{4}{3}|S_{1}|+\frac{1}{3}-c_{1}(G-S_{1})-\frac{2}{3}c_{3}(G-S_{1})=(\frac{4}{3}|S|+\frac{4}{3})+\frac{1}{3}-(c_{1}(G-S)+2)-\frac{2}{3}(c_{3}(G-S)-1)=\beta $.
This contradicts the maximality of $S$.
\qed
\end{enumerate}

Set $T_{1}=\C_{1}(G-S)$, $T_{2}=\C_{3}(G-S)$ and $T=T_{1}\cup T_{2}$.
Now we construct a bipartite graph $H$ with bipartition $(S,T)$ by letting $uC\in E(H)~(u\in S,C\in T)$ if and only if $N_{G}(u)\cap V(C)\not=\emptyset $.

\begin{claim}%%%%%%%%%%%%%%%%%%%%%%%%%%%%%%%%%%%%%%%%%%%%%%%%%%%%%%%%%%%%%%%%%%%%%%%%%%%%%%%%%%%%%%%%%%%%%%%%%%%%%%%%%%%
\label{cl3.2}
The following hold.
\begin{enumerate}[{\upshape(i)}]
\item
$|T_{1}|+\frac{2}{3}|T_{2}|\leq \frac{4}{3}|S|+\frac{1}{3}$.
\item
$1\leq |S|\leq |T_{1}|+|T_{2}|$.
\item
For every $X\subseteq S_{0}$, either $|N_{H}(X)\cap T_{1}|+\frac{2}{3}|N_{H}(X)\cap T_{2}|\geq \frac{4}{3}|X|$ or $N_{H}(X)=T$.
\end{enumerate}
\end{claim}
%%%%%%%%%%%%%%%%%%%%%%%%%%%%%%%%%%%%%%%%%%%%%%%%%%%%%%%%%%%%%%%%%%%%%%%%%%%%%%%%%%%%%%%%%%%%%%%%%%%%%%%%%%%%%%%%%%%%%%%%
\proof
\begin{enumerate}[{\upshape(i)}]
\item
By the assumption of the theorem, $|T_{1}|+\frac{2}{3}|T_{2}|=c_{1}(G-S)+\frac{2}{3}c_{3}(G-S)\leq \frac{4}{3}|S|+\frac{1}{3}$.
\item
Since $c_{1}(G)+c_{3}(G)=0$ and $c_{1}(G-S)+c_{3}(G-S)\geq 1$, $S\not=\emptyset $ (i.e., $|S|\geq 1$).
Since $\frac{4}{3}|S|+\frac{1}{3}-|T_{1}|-\frac{2}{3}|T_{1}|=\frac{4}{3}|S|+\frac{1}{3}-c_{1}(G-S)-\frac{2}{3}c_{3}(G-S)=\beta \leq \frac{5}{3}$, we get $|S|\leq \frac{3}{4}|T_{1}|+\frac{2}{4}|T_{2}|+1\leq \frac{3}{4}|T|+1<|T|+1$, and hence $|S|\leq |T|=|T_{1}|+|T_{2}|$.
\item
Suppose that there exists a set $X\subseteq S$ such that $|N_{H}(X)\cap T_{1}|+\frac{2}{3}|N_{H}(X)\cap T_{2}|<\frac{4}{3}|X|$ and $N_{H}(X)\not=T$.
Since $T-N_{H}(X)\subseteq \C_{1}(G-(S-X))\cup \C_{3}(G-(S-X))$ by the definition of $H$, we have $S-X\in \S$.
We also get $c_{1}(G-(S-X))+\frac{2}{3}c_{3}(G-(S-X))\geq (|T_{1}|-|N_{H}(X)\cap T_{1}|)+\frac{2}{3}(|T_{2}|-|N_{H}(X)\cap T_{2}|)=(c_{1}(G-S)+\frac{2}{3}c_{3}(G-S))-(|N_{H}(X)\cap T_{1}|+\frac{2}{3}|N_{H}(X)\cap T_{2}|)$.
Consequently $\frac{4}{3}|S-X|+\frac{1}{3}-c_{1}(G-(S-X))-\frac{2}{3}c_{3}(G-(S-X))\leq (\frac{4}{3}|S|+\frac{1}{3}-c_{1}(G-S)-\frac{2}{3}c_{3}(G-S))-(\frac{4}{3}|X|-|N_{H}(X)\cap T_{1}|-\frac{2}{3}|N_{H}(X)\cap T_{2}|)<\beta $, which contradicts the definition of $\beta $.
\qed
\end{enumerate}

By Claim~\ref{cl3.2} and Theorem~\ref{thm2}, $H$ has an $S$-central path-factor $F$ such that $V(A)\cap T_{2}\not=\emptyset $ for every $A\in \C_{3}(F)$.
For $A\in \C(F)$, let $U_{A}=V(A)\cap S$, $\L_{A,h}=V(A)\cap T_{h}~(h\in \{1,2\})$, and $\L_{A}=\L_{A,1}\cup \L_{A,2}$.
Let $G_{A}$ be the graph obtained from $G[U_{A}\cup (\bigcup _{C\in \L_{A}}V(C))]$ by deleting all edges of $G[U_{A}]$.

\begin{claim}%%%%%%%%%%%%%%%%%%%%%%%%%%%%%%%%%%%%%%%%%%%%%%%%%%%%%%%%%%%%%%%%%%%%%%%%%%%%%%%%%%%%%%%%%%%%%%%%%%%%%%%%%%%
\label{cl3.3}
For each $A\in \C(F)$, $G_{A}$ has a $\{P_{2},P_{5}\}$-factor.
\end{claim}
%%%%%%%%%%%%%%%%%%%%%%%%%%%%%%%%%%%%%%%%%%%%%%%%%%%%%%%%%%%%%%%%%%%%%%%%%%%%%%%%%%%%%%%%%%%%%%%%%%%%%%%%%%%%%%%%%%%%%%%%
\proof
Since $A$ is a path of $H$, there exists a path $Q_{A}$ of $G_{A}$ such that $U_{A}\subseteq V(Q_{A})$ and $V(Q_{A})\cap V(C)\not=\emptyset $ for every $C\in \L_{A}$.
Choose $Q_{A}$ so that $|V(Q_{A})|$ is as large as possible.
Then for each $C\in \L_{A,2}$ (i.e., $C\in \L_{A}$ with $|V(C)|=3$), since $C$ is complete by Claim~\ref{cl3.1}(ii), it follows that
\begin{align*}
\mbox{either $V(C)\subseteq V(Q_{A})$ or $|V(C)\cap V(Q_{A})|=1$,}
\end{align*}
and
\begin{align*}
\mbox{if $C\in \L_{A}$ is an endvertex of the path $A$ of $H$, then $V(C)\subseteq V(Q_{A})$.}
\end{align*}
Recall that $\L_{A,2}\not=\emptyset $ if $|V(A)|=3$.
Consequently $|V(Q_{A})|\geq |V(A)|$ and, in the case where $|V(A)|\geq 3$, we have $|V(Q_{A})|=5$ or $7$.
Since $|V(A)|\geq 2$, this means that $|V(Q_{A})|\geq 2$ and $|V(Q_{A})|\not=3$.
Furthermore, for each $C\in \L_{A}$, $C-V(Q_{A})$ is either empty or a path of order two.
Therefore if we set $F_{A}=Q_{A}\cup (\bigcup _{C\in \L_{A}}(C-V(Q_{A})))$, then $F_{A}$ is a path-factor of $G_{A}$ with $\C_{3}(F_{A})=\emptyset $.
By Fact~\ref{fact1}, $G_{A}$ has a $\{P_{2},P_{5}\}$-factor.
\qed

By Claims~\ref{cl3.1}(i) and \ref{cl3.3}, $G$ has a $\{P_{2},P_{5}\}$-factor.

This completes the proof of Theorem~\ref{thm1}.

%%%%%%%%%%%%%%%%%%%%%%%%%%%%%%%%%%%%%%%%%%%%%%%%%%%%%%%%%%%%%%%%%%%%%%%%%%%%%%%%%%%%%%%%%%%%%%%%%%%%%%%%%%%%%%%%%%%%%%%%
%%%%%%%%%%%%%%%%%%%%%%%%%%%%%%%%%%%%%%%%%%%%%%%%%%%%%%%%%%%%%%%%%%%%%%%%%%%%%%%%%%%%%%%%%%%%%%%%%%%%%%%%%%%%%%%%%%%%%%%%
\section{Examples}\label{sec4}
%%%%%%%%%%%%%%%%%%%%%%%%%%%%%%%%%%%%%%%%%%%%%%%%%%%%%%%%%%%%%%%%%%%%%%%%%%%%%%%%%%%%%%%%%%%%%%%%%%%%%%%%%%%%%%%%%%%%%%%%
%%%%%%%%%%%%%%%%%%%%%%%%%%%%%%%%%%%%%%%%%%%%%%%%%%%%%%%%%%%%%%%%%%%%%%%%%%%%%%%%%%%%%%%%%%%%%%%%%%%%%%%%%%%%%%%%%%%%%%%%

In this section, we construct graphs having no $\{P_{2},P_{2k+1}\}$-factor.

%%%%%%%%%%%%%%%%%%%%%%%%%%%%%%%%%%%%%%%%%%%%%%%%%%%%%%%%%%%%%%%%%%%%%%%%%%%%%%%%%%%%%%%%%%%%%%%%%%%%%%%%%%%%%%%%%%%%%%%%
%%%%%%%%%%%%%%%%%%%%%%%%%%%%%%%%%%%%%%%%%%%%%%%%%%%%%%%%%%%%%%%%%%%%%%%%%%%%%%%%%%%%%%%%%%%%%%%%%%%%%%%%%%%%%%%%%%%%%%%%
\subsection{Graphs without $\{P_{2},P_{5}\}$-factor}\label{sec4.1}
%%%%%%%%%%%%%%%%%%%%%%%%%%%%%%%%%%%%%%%%%%%%%%%%%%%%%%%%%%%%%%%%%%%%%%%%%%%%%%%%%%%%%%%%%%%%%%%%%%%%%%%%%%%%%%%%%%%%%%%%
%%%%%%%%%%%%%%%%%%%%%%%%%%%%%%%%%%%%%%%%%%%%%%%%%%%%%%%%%%%%%%%%%%%%%%%%%%%%%%%%%%%%%%%%%%%%%%%%%%%%%%%%%%%%%%%%%%%%%%%%

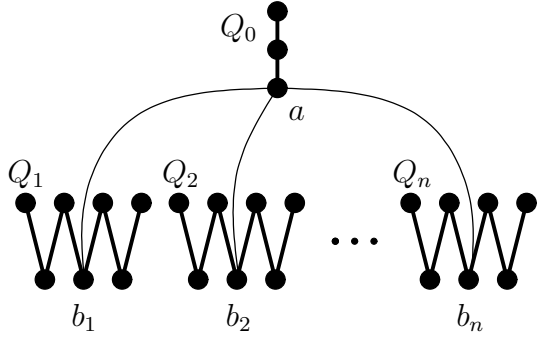
\begin{figure}
\begin{center}
%WinTpicVersion4.26
\unitlength 0.1in
\begin{picture}( 28.5000, 15.8500)(  7.9000,-19.3500)
% CIRCLE 2 0 0 0 Black Black
% 4 995 1400 995 1450 995 1450 995 1450
% 
\special{sh 1.000}%
\special{ia 996 1400 50 50  0.0000000  6.2831853}%
\special{pn 8}%
\special{pn 8}%
\special{ar 996 1400 50 50  0.0000000  6.2831853}%
% CIRCLE 2 0 0 0 Black Black
% 4 1095 1800 1095 1850 1095 1850 1095 1850
% 
\special{sh 1.000}%
\special{ia 1096 1800 50 50  0.0000000  6.2831853}%
\special{pn 8}%
\special{pn 8}%
\special{ar 1096 1800 50 50  0.0000000  6.2831853}%
% CIRCLE 2 0 0 0 Black Black
% 4 1195 1400 1195 1450 1195 1450 1195 1450
% 
\special{sh 1.000}%
\special{ia 1196 1400 50 50  0.0000000  6.2831853}%
\special{pn 8}%
\special{pn 8}%
\special{ar 1196 1400 50 50  0.0000000  6.2831853}%
% CIRCLE 2 0 0 0 Black Black
% 4 1395 1400 1395 1450 1395 1450 1395 1450
% 
\special{sh 1.000}%
\special{ia 1396 1400 50 50  0.0000000  6.2831853}%
\special{pn 8}%
\special{pn 8}%
\special{ar 1396 1400 50 50  0.0000000  6.2831853}%
% CIRCLE 2 0 0 0 Black Black
% 4 1595 1400 1595 1450 1595 1450 1595 1450
% 
\special{sh 1.000}%
\special{ia 1596 1400 50 50  0.0000000  6.2831853}%
\special{pn 8}%
\special{pn 8}%
\special{ar 1596 1400 50 50  0.0000000  6.2831853}%
% CIRCLE 2 0 0 0 Black Black
% 4 1295 1800 1295 1850 1295 1850 1295 1850
% 
\special{sh 1.000}%
\special{ia 1296 1800 50 50  0.0000000  6.2831853}%
\special{pn 8}%
\special{pn 8}%
\special{ar 1296 1800 50 50  0.0000000  6.2831853}%
% CIRCLE 2 0 0 0 Black Black
% 4 1495 1800 1495 1850 1495 1850 1495 1850
% 
\special{sh 1.000}%
\special{ia 1496 1800 50 50  0.0000000  6.2831853}%
\special{pn 8}%
\special{pn 8}%
\special{ar 1496 1800 50 50  0.0000000  6.2831853}%
% LINE 0 0 3 0 Black Black
% 2 995 1400 1095 1800
% 
\special{pn 20}%
\special{pa 996 1400}%
\special{pa 1096 1800}%
\special{fp}%
% LINE 0 0 3 0 Black Black
% 2 1195 1400 1295 1800
% 
\special{pn 20}%
\special{pa 1196 1400}%
\special{pa 1296 1800}%
\special{fp}%
% LINE 0 0 3 0 Black Black
% 2 1395 1400 1495 1800
% 
\special{pn 20}%
\special{pa 1396 1400}%
\special{pa 1496 1800}%
\special{fp}%
% LINE 0 0 3 0 Black Black
% 2 1595 1400 1495 1800
% 
\special{pn 20}%
\special{pa 1596 1400}%
\special{pa 1496 1800}%
\special{fp}%
% LINE 0 0 3 0 Black Black
% 2 1395 1400 1295 1800
% 
\special{pn 20}%
\special{pa 1396 1400}%
\special{pa 1296 1800}%
\special{fp}%
% LINE 0 0 3 0 Black Black
% 2 1195 1400 1095 1800
% 
\special{pn 20}%
\special{pa 1196 1400}%
\special{pa 1096 1800}%
\special{fp}%
% CIRCLE 2 0 0 0 Black Black
% 4 1790 1400 1790 1450 1790 1450 1790 1450
% 
\special{sh 1.000}%
\special{ia 1790 1400 50 50  0.0000000  6.2831853}%
\special{pn 8}%
\special{pn 8}%
\special{ar 1790 1400 50 50  0.0000000  6.2831853}%
% CIRCLE 2 0 0 0 Black Black
% 4 1890 1800 1890 1850 1890 1850 1890 1850
% 
\special{sh 1.000}%
\special{ia 1890 1800 50 50  0.0000000  6.2831853}%
\special{pn 8}%
\special{pn 8}%
\special{ar 1890 1800 50 50  0.0000000  6.2831853}%
% CIRCLE 2 0 0 0 Black Black
% 4 1990 1400 1990 1450 1990 1450 1990 1450
% 
\special{sh 1.000}%
\special{ia 1990 1400 50 50  0.0000000  6.2831853}%
\special{pn 8}%
\special{pn 8}%
\special{ar 1990 1400 50 50  0.0000000  6.2831853}%
% CIRCLE 2 0 0 0 Black Black
% 4 2190 1400 2190 1450 2190 1450 2190 1450
% 
\special{sh 1.000}%
\special{ia 2190 1400 50 50  0.0000000  6.2831853}%
\special{pn 8}%
\special{pn 8}%
\special{ar 2190 1400 50 50  0.0000000  6.2831853}%
% CIRCLE 2 0 0 0 Black Black
% 4 2390 1400 2390 1450 2390 1450 2390 1450
% 
\special{sh 1.000}%
\special{ia 2390 1400 50 50  0.0000000  6.2831853}%
\special{pn 8}%
\special{pn 8}%
\special{ar 2390 1400 50 50  0.0000000  6.2831853}%
% CIRCLE 2 0 0 0 Black Black
% 4 2090 1800 2090 1850 2090 1850 2090 1850
% 
\special{sh 1.000}%
\special{ia 2090 1800 50 50  0.0000000  6.2831853}%
\special{pn 8}%
\special{pn 8}%
\special{ar 2090 1800 50 50  0.0000000  6.2831853}%
% CIRCLE 2 0 0 0 Black Black
% 4 2290 1800 2290 1850 2290 1850 2290 1850
% 
\special{sh 1.000}%
\special{ia 2290 1800 50 50  0.0000000  6.2831853}%
\special{pn 8}%
\special{pn 8}%
\special{ar 2290 1800 50 50  0.0000000  6.2831853}%
% LINE 0 0 3 0 Black Black
% 2 1790 1400 1890 1800
% 
\special{pn 20}%
\special{pa 1790 1400}%
\special{pa 1890 1800}%
\special{fp}%
% LINE 0 0 3 0 Black Black
% 2 1990 1400 2090 1800
% 
\special{pn 20}%
\special{pa 1990 1400}%
\special{pa 2090 1800}%
\special{fp}%
% LINE 0 0 3 0 Black Black
% 2 2190 1400 2290 1800
% 
\special{pn 20}%
\special{pa 2190 1400}%
\special{pa 2290 1800}%
\special{fp}%
% LINE 0 0 3 0 Black Black
% 2 2390 1400 2290 1800
% 
\special{pn 20}%
\special{pa 2390 1400}%
\special{pa 2290 1800}%
\special{fp}%
% LINE 0 0 3 0 Black Black
% 2 2190 1400 2090 1800
% 
\special{pn 20}%
\special{pa 2190 1400}%
\special{pa 2090 1800}%
\special{fp}%
% LINE 0 0 3 0 Black Black
% 2 1990 1400 1890 1800
% 
\special{pn 20}%
\special{pa 1990 1400}%
\special{pa 1890 1800}%
\special{fp}%
% CIRCLE 2 0 0 0 Black Black
% 4 2990 1400 2990 1450 2990 1450 2990 1450
% 
\special{sh 1.000}%
\special{ia 2990 1400 50 50  0.0000000  6.2831853}%
\special{pn 8}%
\special{pn 8}%
\special{ar 2990 1400 50 50  0.0000000  6.2831853}%
% CIRCLE 2 0 0 0 Black Black
% 4 3090 1800 3090 1850 3090 1850 3090 1850
% 
\special{sh 1.000}%
\special{ia 3090 1800 50 50  0.0000000  6.2831853}%
\special{pn 8}%
\special{pn 8}%
\special{ar 3090 1800 50 50  0.0000000  6.2831853}%
% CIRCLE 2 0 0 0 Black Black
% 4 3190 1400 3190 1450 3190 1450 3190 1450
% 
\special{sh 1.000}%
\special{ia 3190 1400 50 50  0.0000000  6.2831853}%
\special{pn 8}%
\special{pn 8}%
\special{ar 3190 1400 50 50  0.0000000  6.2831853}%
% CIRCLE 2 0 0 0 Black Black
% 4 3390 1400 3390 1450 3390 1450 3390 1450
% 
\special{sh 1.000}%
\special{ia 3390 1400 50 50  0.0000000  6.2831853}%
\special{pn 8}%
\special{pn 8}%
\special{ar 3390 1400 50 50  0.0000000  6.2831853}%
% CIRCLE 2 0 0 0 Black Black
% 4 3590 1400 3590 1450 3590 1450 3590 1450
% 
\special{sh 1.000}%
\special{ia 3590 1400 50 50  0.0000000  6.2831853}%
\special{pn 8}%
\special{pn 8}%
\special{ar 3590 1400 50 50  0.0000000  6.2831853}%
% CIRCLE 2 0 0 0 Black Black
% 4 3290 1800 3290 1850 3290 1850 3290 1850
% 
\special{sh 1.000}%
\special{ia 3290 1800 50 50  0.0000000  6.2831853}%
\special{pn 8}%
\special{pn 8}%
\special{ar 3290 1800 50 50  0.0000000  6.2831853}%
% CIRCLE 2 0 0 0 Black Black
% 4 3490 1800 3490 1850 3490 1850 3490 1850
% 
\special{sh 1.000}%
\special{ia 3490 1800 50 50  0.0000000  6.2831853}%
\special{pn 8}%
\special{pn 8}%
\special{ar 3490 1800 50 50  0.0000000  6.2831853}%
% LINE 0 0 3 0 Black Black
% 2 2990 1400 3090 1800
% 
\special{pn 20}%
\special{pa 2990 1400}%
\special{pa 3090 1800}%
\special{fp}%
% LINE 0 0 3 0 Black Black
% 2 3190 1400 3290 1800
% 
\special{pn 20}%
\special{pa 3190 1400}%
\special{pa 3290 1800}%
\special{fp}%
% LINE 0 0 3 0 Black Black
% 2 3390 1400 3490 1800
% 
\special{pn 20}%
\special{pa 3390 1400}%
\special{pa 3490 1800}%
\special{fp}%
% LINE 0 0 3 0 Black Black
% 2 3590 1400 3490 1800
% 
\special{pn 20}%
\special{pa 3590 1400}%
\special{pa 3490 1800}%
\special{fp}%
% LINE 0 0 3 0 Black Black
% 2 3390 1400 3290 1800
% 
\special{pn 20}%
\special{pa 3390 1400}%
\special{pa 3290 1800}%
\special{fp}%
% LINE 0 0 3 0 Black Black
% 2 3190 1400 3090 1800
% 
\special{pn 20}%
\special{pa 3190 1400}%
\special{pa 3090 1800}%
\special{fp}%
% CIRCLE 2 0 0 0 Black Black
% 4 2300 800 2300 850 2300 850 2300 850
% 
\special{sh 1.000}%
\special{ia 2300 800 50 50  0.0000000  6.2831853}%
\special{pn 8}%
\special{pn 8}%
\special{ar 2300 800 50 50  0.0000000  6.2831853}%
% SPLINE 2 0 3 0 Black Black
% 5 2300 800 1500 1000 1300 1400 1300 1800 1300 1800
% 
\special{pn 8}%
\special{pa 2300 800}%
\special{pa 2266 802}%
\special{pa 2230 802}%
\special{pa 2194 804}%
\special{pa 2160 806}%
\special{pa 2124 806}%
\special{pa 2090 808}%
\special{pa 2054 810}%
\special{pa 2020 814}%
\special{pa 1986 816}%
\special{pa 1952 820}%
\special{pa 1920 824}%
\special{pa 1886 830}%
\special{pa 1822 842}%
\special{pa 1792 850}%
\special{pa 1760 858}%
\special{pa 1730 868}%
\special{pa 1702 878}%
\special{pa 1672 890}%
\special{pa 1644 902}%
\special{pa 1592 930}%
\special{pa 1566 948}%
\special{pa 1518 984}%
\special{pa 1474 1028}%
\special{pa 1454 1052}%
\special{pa 1434 1078}%
\special{pa 1398 1130}%
\special{pa 1382 1160}%
\special{pa 1368 1188}%
\special{pa 1354 1218}%
\special{pa 1342 1250}%
\special{pa 1330 1280}%
\special{pa 1320 1312}%
\special{pa 1312 1344}%
\special{pa 1294 1440}%
\special{pa 1292 1472}%
\special{pa 1288 1504}%
\special{pa 1288 1536}%
\special{pa 1286 1570}%
\special{pa 1288 1602}%
\special{pa 1288 1632}%
\special{pa 1294 1728}%
\special{pa 1298 1760}%
\special{pa 1300 1792}%
\special{pa 1300 1800}%
\special{sp}%
% SPLINE 2 0 3 0 Black Black
% 4 2300 800 2100 1200 2100 1800 2100 1800
% 
\special{pn 8}%
\special{pa 2300 800}%
\special{pa 2284 828}%
\special{pa 2248 884}%
\special{pa 2200 968}%
\special{pa 2184 998}%
\special{pa 2156 1054}%
\special{pa 2142 1084}%
\special{pa 2130 1112}%
\special{pa 2100 1202}%
\special{pa 2092 1232}%
\special{pa 2086 1262}%
\special{pa 2082 1292}%
\special{pa 2074 1356}%
\special{pa 2072 1386}%
\special{pa 2072 1418}%
\special{pa 2070 1450}%
\special{pa 2072 1482}%
\special{pa 2072 1514}%
\special{pa 2074 1548}%
\special{pa 2076 1580}%
\special{pa 2080 1612}%
\special{pa 2082 1646}%
\special{pa 2090 1710}%
\special{pa 2094 1744}%
\special{pa 2098 1776}%
\special{pa 2100 1800}%
\special{sp}%
% SPLINE 2 0 3 0 Black Black
% 5 2300 800 3100 1000 3300 1400 3300 1800 3300 1800
% 
\special{pn 8}%
\special{pa 2300 800}%
\special{pa 2336 802}%
\special{pa 2372 802}%
\special{pa 2406 804}%
\special{pa 2442 806}%
\special{pa 2476 806}%
\special{pa 2512 808}%
\special{pa 2546 810}%
\special{pa 2580 814}%
\special{pa 2614 816}%
\special{pa 2682 824}%
\special{pa 2778 842}%
\special{pa 2810 850}%
\special{pa 2840 858}%
\special{pa 2900 878}%
\special{pa 2956 902}%
\special{pa 2984 916}%
\special{pa 3010 930}%
\special{pa 3036 948}%
\special{pa 3084 984}%
\special{pa 3128 1028}%
\special{pa 3148 1052}%
\special{pa 3166 1078}%
\special{pa 3186 1104}%
\special{pa 3202 1130}%
\special{pa 3218 1160}%
\special{pa 3232 1188}%
\special{pa 3246 1218}%
\special{pa 3258 1250}%
\special{pa 3270 1280}%
\special{pa 3280 1312}%
\special{pa 3296 1376}%
\special{pa 3302 1408}%
\special{pa 3310 1472}%
\special{pa 3314 1536}%
\special{pa 3314 1602}%
\special{pa 3312 1632}%
\special{pa 3312 1664}%
\special{pa 3310 1696}%
\special{pa 3306 1728}%
\special{pa 3302 1792}%
\special{pa 3300 1800}%
\special{sp}%
% STR 2 0 3 0 Black Black
% 4 1000 1150 1000 1250 5 0 0 0
% $Q_{1}$
\put(10.0000,-12.5000){\makebox(0,0){$Q_{1}$}}%
% CIRCLE 2 0 0 0 Black Black
% 4 2300 600 2300 650 2300 650 2300 650
% 
\special{sh 1.000}%
\special{ia 2300 600 50 50  0.0000000  6.2831853}%
\special{pn 8}%
\special{pn 8}%
\special{ar 2300 600 50 50  0.0000000  6.2831853}%
% CIRCLE 2 0 0 0 Black Black
% 4 2300 400 2300 450 2300 450 2300 450
% 
\special{sh 1.000}%
\special{ia 2300 400 50 50  0.0000000  6.2831853}%
\special{pn 8}%
\special{pn 8}%
\special{ar 2300 400 50 50  0.0000000  6.2831853}%
% LINE 0 0 3 0 Black Black
% 2 2300 400 2300 800
% 
\special{pn 20}%
\special{pa 2300 400}%
\special{pa 2300 800}%
\special{fp}%
% STR 2 0 3 0 Black Black
% 4 1300 1900 1300 2000 5 0 0 0
% $b_{1}$
\put(13.0000,-20.0000){\makebox(0,0){$b_{1}$}}%
% STR 2 0 3 0 Black Black
% 4 2100 390 2100 490 5 0 0 0
% $Q_{0}$
\put(21.0000,-4.9000){\makebox(0,0){$Q_{0}$}}%
% STR 2 0 3 0 Black Black
% 4 2400 830 2400 930 5 0 0 0
% $a$
\put(24.0000,-9.3000){\makebox(0,0){$a$}}%
% DOT 0 0 3 0 Black Black
% 4 2800 1600 2700 1600 2600 1600 2600 1600
% 
\special{pn 4}%
\special{sh 1}%
\special{ar 2800 1600 16 16 0  6.28318530717959E+0000}%
\special{sh 1}%
\special{ar 2700 1600 16 16 0  6.28318530717959E+0000}%
\special{sh 1}%
\special{ar 2600 1600 16 16 0  6.28318530717959E+0000}%
\special{sh 1}%
\special{ar 2600 1600 16 16 0  6.28318530717959E+0000}%
% STR 2 0 3 0 Black Black
% 4 3300 1900 3300 2000 5 0 0 0
% $b_{n}$
\put(33.0000,-20.0000){\makebox(0,0){$b_{n}$}}%
% STR 2 0 3 0 Black Black
% 4 1800 1150 1800 1250 5 0 0 0
% $Q_{2}$
\put(18.0000,-12.5000){\makebox(0,0){$Q_{2}$}}%
% STR 2 0 3 0 Black Black
% 4 3000 1150 3000 1250 5 0 0 0
% $Q_{n}$
\put(30.0000,-12.5000){\makebox(0,0){$Q_{n}$}}%
% STR 2 0 3 0 Black Black
% 4 2100 1900 2100 2000 5 0 0 0
% $b_{2}$
\put(21.0000,-20.0000){\makebox(0,0){$b_{2}$}}%
\end{picture}%
\caption{Graph $H_{n}$}
\label{thm3.1fig1}
\end{center}
\end{figure}

Let $n\geq 1$ be an integer.
Let $Q_{0}$ be a path of order $3$, and let $a$ be an endvertex of $Q_{0}$.
Let $Q_{1},\ldots ,Q_{n}$ be disjoint paths of order $7$, and for each $i~(1\leq i\leq n)$, let $b_{i}$ be the center of $Q_{i}$.
Let $H_{n}$ denote the graph obtained from $\bigcup _{0\leq i\leq n}Q_{i}$ by joining $a$ to $b_{i}$ for every $i~(1\leq i\leq n)$ (see Figure~\ref{thm3.1fig1}).

Suppose that $H_{n}$ has a $\{P_{2},P_{5}\}$-factor $F$.
Since $Q_{0}$ does not have a $\{P_{2},P_{5}\}$-factor, $F$ contains $ab_{i}$ for some $i~(1\leq i\leq n)$.
Since $d_{F}(b_{i})\leq 2$, this requires that at least one of the components of $Q_{i}-b_{i}$ should have a $\{P_{2},P_{5}\}$-factor, which is impossible because each component of $Q_{i}-b_{i}$ is a path of order $3$.
Thus $H_{n}$ has no $\{P_{2},P_{5}\}$-factor.

\begin{lem}%%%%%%%%%%%%%%%%%%%%%%%%%%%%%%%%%%%%%%%%%%%%%%%%%%%%%%%%%%%%%%%%%%%%%%%%%%%%%%%%%%%%%%%%%%%%%%%%%%%%%%%%%%%%%
\label{lem4.1.1}
For all $X\subseteq V(H_{n})$, $c_{1}(H_{n}-X)+\frac{2}{3}c_{3}(H_{n}-X)\leq \frac{4}{3}|X|+\frac{2}{3}$.
\end{lem}
%%%%%%%%%%%%%%%%%%%%%%%%%%%%%%%%%%%%%%%%%%%%%%%%%%%%%%%%%%%%%%%%%%%%%%%%%%%%%%%%%%%%%%%%%%%%%%%%%%%%%%%%%%%%%%%%%%%%%%%%
\proof
Let $X\subseteq V(H_{n})$.
Then we can verify that
\begin{align}
c_{1}(Q_{0}-X)+\frac{2}{3}c_{3}(Q_{0}-X)\leq \frac{4}{3}|V(Q_{0})\cap X|+\frac{2}{3}\label{sec4-1}
\end{align}
and
\begin{align}
c_{1}(Q_{i}-X)+\frac{2}{3}c_{3}(Q_{i}-X)\leq \frac{4}{3}|V(Q_{i})\cap X|\mbox{ for every }i~(1\leq i\leq n)\label{sec3-2}
\end{align}
Since every component $C$ of $H_{n}-X$ with $|V(C)|=1$ belongs to $\bigcup _{0\leq i\leq n}\C_{1}(Q_{i}-X)$, we have
\begin{align}
|\C_{1}(H_{n}-X)|=\sum _{0\leq i\leq n}|\C_{1}(Q_{i}-X)|-\left|\left(\bigcup _{0\leq i\leq n}\C_{1}(Q_{i}-X)\right)-\C_{1}(H_{n}-X)\right|.\label{sec4-3}
\end{align}
Furthermore,
\begin{align}
|\C_{3}(H_{n}-X)|\leq \sum _{0\leq i\leq n}|\C_{3}(Q_{i}-X)|+\left|\C_{3}(H_{n}-X)-\left(\bigcup _{0\leq i\leq n}\C_{3}(Q_{i}-X)\right)\right|.\label{sec4-4}
\end{align}

Let $C$ be a component of $H_{n}-X$ with $|V(C)|=3$ which does not belong to $\bigcup _{0\leq i\leq n}\C_{3}(Q_{i}-X)$.
Then $C$ intersects with at least two of the $Q_{i}~(0\leq i\leq n)$.
Since $|V(C)|=3$, $C$ contains a component of $Q_{i}-X$ of order $1$ for some $i~(0\leq i\leq n)$.
Since $C$ is arbitrary, this implies that
\begin{align}
\left|\C_{3}(H_{n}-X)-\left(\bigcup _{0\leq i\leq n}\C_{3}(Q_{i}-X)\right)\right|\leq \left|\left(\bigcup _{0\leq i\leq n}\C_{1}(Q_{i}-X)\right)-\C_{1}(H_{n}-X)\right|.\label{sec4-5}
\end{align}
By (\ref{sec4-1})--(\ref{sec4-5}),
\begin{align}
c_{1}(H_{n}-X)&+\frac{2}{3}c_{3}(H_{n}-X)\nonumber \\
&\leq \left(\sum _{0\leq i\leq n}|\C_{1}(Q_{i}-X)|-\left|\left(\bigcup _{0\leq i\leq n}\C_{1}(Q_{i}-X)\right)-\C_{1}(H_{n}-X)\right|\right)\nonumber \\
&\quad +\frac{2}{3}\left(\sum _{0\leq i\leq n}|\C_{3}(Q_{i}-X)|+\left|\C_{3}(H_{n}-X)-\left(\bigcup _{0\leq i\leq n}\C_{3}(Q_{i}-X)\right)\right|\right)\nonumber \\
&\leq \left(\sum _{0\leq i\leq n}|\C_{1}(Q_{i}-X)|-\left|\left(\bigcup _{0\leq i\leq n}\C_{1}(Q_{i}-X)\right)-\C_{1}(H_{n}-X)\right|\right)\nonumber \\
&\quad +\frac{2}{3}\left(\sum _{0\leq i\leq n}|\C_{3}(Q_{i}-X)|+\left|\left(\bigcup _{0\leq i\leq n}\C_{1}(Q_{i}-X)\right)-\C_{1}(H_{n}-X)\right|\right)\nonumber \\
&\leq \sum _{0\leq i\leq n}|\C_{1}(Q_{i}-X)|+\frac{2}{3}\sum _{0\leq i\leq n}|\C_{3}(Q_{i}-X)|\nonumber \\
&= \sum _{0\leq i\leq n}\left(c_{1}(Q_{i}-X)+\frac{2}{3}c_{3}(Q_{i}-X)\right)\nonumber \\
&\leq \frac{4}{3}\sum _{0\leq i\leq n}|V(Q_{i})\cap X|+\frac{2}{3}\nonumber \\
&= \frac{4}{3}|X|+\frac{2}{3}.\nonumber
\end{align}
Thus we get the desired conclusion.
\qed

From Lemma~\ref{lem4.1.1}, we get the following proposition, which implies that Theorem~\ref{thm1} is best possible.

\begin{prop}%%%%%%%%%%%%%%%%%%%%%%%%%%%%%%%%%%%%%%%%%%%%%%%%%%%%%%%%%%%%%%%%%%%%%%%%%%%%%%%%%%%%%%%%%%%%%%%%%%%%%%%%%%%%
\label{prop3.1}
There exist infinitely many graphs $G$ having no $\{P_{2},P_{5}\}$-factor such that $c_{1}(G-X)+\frac{2}{3}c_{3}(G-X)\leq \frac{4}{3}|X|+\frac{2}{3}$ for all $X\subseteq V(G)$.
\end{prop}
%%%%%%%%%%%%%%%%%%%%%%%%%%%%%%%%%%%%%%%%%%%%%%%%%%%%%%%%%%%%%%%%%%%%%%%%%%%%%%%%%%%%%%%%%%%%%%%%%%%%%%%%%%%%%%%%%%%%%%%%

%%%%%%%%%%%%%%%%%%%%%%%%%%%%%%%%%%%%%%%%%%%%%%%%%%%%%%%%%%%%%%%%%%%%%%%%%%%%%%%%%%%%%%%%%%%%%%%%%%%%%%%%%%%%%%%%%%%%%%%%
%%%%%%%%%%%%%%%%%%%%%%%%%%%%%%%%%%%%%%%%%%%%%%%%%%%%%%%%%%%%%%%%%%%%%%%%%%%%%%%%%%%%%%%%%%%%%%%%%%%%%%%%%%%%%%%%%%%%%%%%
\subsection{Graphs without $\{P_{2},P_{2k+1}\}$-factor for $k\geq 3$}\label{sec4.2}
%%%%%%%%%%%%%%%%%%%%%%%%%%%%%%%%%%%%%%%%%%%%%%%%%%%%%%%%%%%%%%%%%%%%%%%%%%%%%%%%%%%%%%%%%%%%%%%%%%%%%%%%%%%%%%%%%%%%%%%%
%%%%%%%%%%%%%%%%%%%%%%%%%%%%%%%%%%%%%%%%%%%%%%%%%%%%%%%%%%%%%%%%%%%%%%%%%%%%%%%%%%%%%%%%%%%%%%%%%%%%%%%%%%%%%%%%%%%%%%%%

\begin{figure}
\begin{center}
%WinTpicVersion4.28b
{\unitlength 0.1in
\begin{picture}( 31.0000, 15.3500)( 12.5000,-21.8500)
% CIRCLE 2 0 0 0 Black Black
% 4 1400 2000 1400 2050 1400 2050 1400 2050
% 
\special{sh 1.000}%
\special{ia 1400 2000 50 50  0.0000000  6.2831853}%
\special{pn 8}%
\special{ar 1400 2000 50 50  0.0000000  6.2831853}%
% CIRCLE 2 0 0 0 Black Black
% 4 1600 2000 1600 2050 1600 2050 1600 2050
% 
\special{sh 1.000}%
\special{ia 1600 2000 50 50  0.0000000  6.2831853}%
\special{pn 8}%
\special{ar 1600 2000 50 50  0.0000000  6.2831853}%
% CIRCLE 2 0 0 0 Black Black
% 4 2200 2000 2200 2050 2200 2050 2200 2050
% 
\special{sh 1.000}%
\special{ia 2200 2000 50 50  0.0000000  6.2831853}%
\special{pn 8}%
\special{ar 2200 2000 50 50  0.0000000  6.2831853}%
% CIRCLE 2 0 0 0 Black Black
% 4 2400 2000 2400 2050 2400 2050 2400 2050
% 
\special{sh 1.000}%
\special{ia 2400 2000 50 50  0.0000000  6.2831853}%
\special{pn 8}%
\special{ar 2400 2000 50 50  0.0000000  6.2831853}%
% DOT 0 0 3 0 Black Black
% 4 2000 2000 1900 2000 1800 2000 1800 2000
% 
\special{pn 4}%
\special{sh 1}%
\special{ar 2000 2000 16 16 0  6.28318530717959E+0000}%
\special{sh 1}%
\special{ar 1900 2000 16 16 0  6.28318530717959E+0000}%
\special{sh 1}%
\special{ar 1800 2000 16 16 0  6.28318530717959E+0000}%
\special{sh 1}%
\special{ar 1800 2000 16 16 0  6.28318530717959E+0000}%
% LINE 1 0 3 0 Black Black
% 4 2200 2000 2400 2000 1600 2000 1400 2000
% 
\special{pn 13}%
\special{pa 2200 2000}%
\special{pa 2400 2000}%
\special{fp}%
\special{pa 1600 2000}%
\special{pa 1400 2000}%
\special{fp}%
% STR 2 0 3 0 Black Black
% 4 1900 1700 1900 1800 5 0 0 0
% $+$
\put(19.0000,-18.0000){\makebox(0,0){$+$}}%
% CIRCLE 2 0 0 0 Black Black
% 4 3200 2000 3200 2050 3200 2050 3200 2050
% 
\special{sh 1.000}%
\special{ia 3200 2000 50 50  0.0000000  6.2831853}%
\special{pn 8}%
\special{ar 3200 2000 50 50  0.0000000  6.2831853}%
% CIRCLE 2 0 0 0 Black Black
% 4 3400 2000 3400 2050 3400 2050 3400 2050
% 
\special{sh 1.000}%
\special{ia 3400 2000 50 50  0.0000000  6.2831853}%
\special{pn 8}%
\special{ar 3400 2000 50 50  0.0000000  6.2831853}%
% CIRCLE 2 0 0 0 Black Black
% 4 4000 2000 4000 2050 4000 2050 4000 2050
% 
\special{sh 1.000}%
\special{ia 4000 2000 50 50  0.0000000  6.2831853}%
\special{pn 8}%
\special{ar 4000 2000 50 50  0.0000000  6.2831853}%
% CIRCLE 2 0 0 0 Black Black
% 4 4200 2000 4200 2050 4200 2050 4200 2050
% 
\special{sh 1.000}%
\special{ia 4200 2000 50 50  0.0000000  6.2831853}%
\special{pn 8}%
\special{ar 4200 2000 50 50  0.0000000  6.2831853}%
% DOT 0 0 3 0 Black Black
% 4 3800 2000 3700 2000 3600 2000 3600 2000
% 
\special{pn 4}%
\special{sh 1}%
\special{ar 3800 2000 16 16 0  6.28318530717959E+0000}%
\special{sh 1}%
\special{ar 3700 2000 16 16 0  6.28318530717959E+0000}%
\special{sh 1}%
\special{ar 3600 2000 16 16 0  6.28318530717959E+0000}%
\special{sh 1}%
\special{ar 3600 2000 16 16 0  6.28318530717959E+0000}%
% LINE 0 0 3 0 Black Black
% 4 4000 2000 4200 2000 3400 2000 3200 2000
% 
\special{pn 20}%
\special{pa 4000 2000}%
\special{pa 4200 2000}%
\special{fp}%
\special{pa 3400 2000}%
\special{pa 3200 2000}%
\special{fp}%
% STR 2 0 3 0 Black Black
% 4 3700 1700 3700 1800 5 0 0 0
% $+$
\put(37.0000,-18.0000){\makebox(0,0){$+$}}%
% DOT 0 0 3 0 Black Black
% 4 2900 1800 2800 1800 2700 1800 2700 1800
% 
\special{pn 4}%
\special{sh 1}%
\special{ar 2900 1800 16 16 0  6.28318530717959E+0000}%
\special{sh 1}%
\special{ar 2800 1800 16 16 0  6.28318530717959E+0000}%
\special{sh 1}%
\special{ar 2700 1800 16 16 0  6.28318530717959E+0000}%
\special{sh 1}%
\special{ar 2700 1800 16 16 0  6.28318530717959E+0000}%
% STR 2 0 3 0 Black Black
% 4 3700 2150 3700 2250 5 0 0 0
% $R_{2n+1}$
\put(37.0000,-22.5000){\makebox(0,0){$R_{2n+1}$}}%
% BOX 2 0 3 0 Black Black
% 2 1250 1450 2550 2150
% 
\special{pn 8}%
\special{pa 1250 1450}%
\special{pa 2550 1450}%
\special{pa 2550 2150}%
\special{pa 1250 2150}%
\special{pa 1250 1450}%
\special{pa 2550 1450}%
\special{fp}%
% BOX 2 0 3 0 Black Black
% 2 1300 1900 2500 2100
% 
\special{pn 8}%
\special{pa 1300 1900}%
\special{pa 2500 1900}%
\special{pa 2500 2100}%
\special{pa 1300 2100}%
\special{pa 1300 1900}%
\special{pa 2500 1900}%
\special{fp}%
% BOX 2 0 3 0 Black Black
% 2 3100 1900 4300 2100
% 
\special{pn 8}%
\special{pa 3100 1900}%
\special{pa 4300 1900}%
\special{pa 4300 2100}%
\special{pa 3100 2100}%
\special{pa 3100 1900}%
\special{pa 4300 1900}%
\special{fp}%
% BOX 2 0 3 0 Black Black
% 2 3050 1450 4350 2150
% 
\special{pn 8}%
\special{pa 3050 1450}%
\special{pa 4350 1450}%
\special{pa 4350 2150}%
\special{pa 3050 2150}%
\special{pa 3050 1450}%
\special{pa 4350 1450}%
\special{fp}%
% STR 2 0 3 0 Black Black
% 4 1900 2150 1900 2250 5 0 0 0
% $R_{1}$
\put(19.0000,-22.5000){\makebox(0,0){$R_{1}$}}%
% ELLIPSE 2 0 3 0 Black Black
% 4 2800 850 3200 1050 3200 1050 3200 1050
% 
\special{pn 8}%
\special{ar 2800 850 400 200  0.0000000  6.2831853}%
% STR 2 0 3 0 Black Black
% 4 2800 750 2800 850 5 0 0 0
% $R_{0}$
\put(28.0000,-8.5000){\makebox(0,0){$R_{0}$}}%
% ELLIPSE 2 0 3 0 Black Black
% 4 1900 1600 2400 1700 2400 1700 2400 1700
% 
\special{pn 8}%
\special{ar 1900 1600 500 100  0.0000000  6.2831853}%
% STR 2 0 3 0 Black Black
% 4 1900 1500 1900 1600 5 0 0 0
% $K_{1}$
\put(19.0000,-16.0000){\makebox(0,0){$K_{1}$}}%
% ELLIPSE 2 0 3 0 Black Black
% 4 3700 1600 4200 1700 4200 1700 4200 1700
% 
\special{pn 8}%
\special{ar 3700 1600 500 100  0.0000000  6.2831853}%
% STR 2 0 3 0 Black Black
% 4 3700 1500 3700 1600 5 0 0 0
% $K_{2n+1}$
\put(37.0000,-16.0000){\makebox(0,0){$K_{2n+1}$}}%
% LINE 2 0 3 0 Black Black
% 2 1250 1450 2430 770
% 
\special{pn 8}%
\special{pa 1250 1450}%
\special{pa 2430 770}%
\special{fp}%
% LINE 2 0 3 0 Black Black
% 2 3170 770 4350 1450
% 
\special{pn 8}%
\special{pa 3170 770}%
\special{pa 4350 1450}%
\special{fp}%
% LINE 2 0 3 0 Black Black
% 2 2550 1450 3160 940
% 
\special{pn 8}%
\special{pa 2550 1450}%
\special{pa 3160 940}%
\special{fp}%
% LINE 2 0 3 0 Black Black
% 2 2440 940 3050 1450
% 
\special{pn 8}%
\special{pa 2440 940}%
\special{pa 3050 1450}%
\special{fp}%
% STR 2 0 3 0 Black Black
% 4 2200 1100 2200 1200 5 0 0 0
% $+$
\put(22.0000,-12.0000){\makebox(0,0){$+$}}%
% STR 2 0 3 0 Black Black
% 4 3400 1100 3400 1200 5 0 0 0
% $+$
\put(34.0000,-12.0000){\makebox(0,0){$+$}}%
\end{picture}}%
\caption{Graph $H'_{n}$}
\label{thm3.2fig1}
\end{center}
\end{figure}
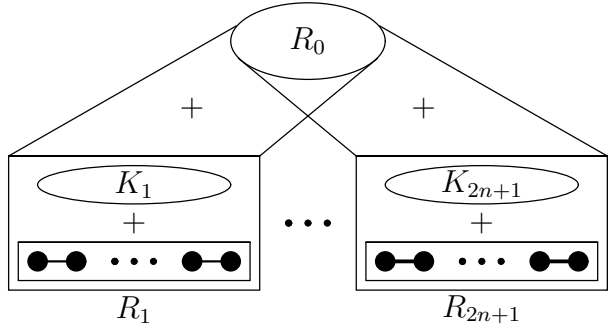

Let $k\geq 3$ be an integer with $k\equiv 0~(\mbox{mod }3)$, and write $k=3m$.
Let $n\geq 1$ be an integer.
Let $R_{0}$ be a complete graph of order $n$.
For each $i~(1\leq i\leq 2n+1)$, let $K_{i}$ be a complete graph of order $2m-1$, and let $R_{i}$ denote the graph obtained from $K_{i}$ by joining each vertex of the union of $2m+1$ disjoint paths of order $2$ to all vertices of $K_{i}$.
Let $H'_{n}=R_{0}+(\bigcup _{1\leq i\leq 2n+1}R_{i})$ (see Figure~\ref{thm3.2fig1}).

Since $|V(R_{i})|=2k+1$ and $R_{i}$ does not contain a path of order $2k+1$, $R_{i}$ has no $\{P_{2},P_{2k+1}\}$-factor.
Suppose that $H'_{n}$ has a $\{P_{2},P_{2k+1}\}$-factor $F$.
Then for each $i~(1\leq i\leq 2n+1)$, $F$ contains an edge joining $V(R_{i})$ and $V(R_{0})$.
Since $2n+1>2|V(R_{0})|$, this implies that there exists $x\in V(R_{0})$ such that $d_{F}(x)\geq 3$, which is a contradiction.
Thus $H'_{n}$ has no $\{P_{2},P_{2k+1}\}$-factor.

\begin{lem}%%%%%%%%%%%%%%%%%%%%%%%%%%%%%%%%%%%%%%%%%%%%%%%%%%%%%%%%%%%%%%%%%%%%%%%%%%%%%%%%%%%%%%%%%%%%%%%%%%%%%%%%%%%%%
\label{lem4.2.1}
For all $X\subseteq V(H'_{n})$, $\sum _{0\leq j\leq k-1}c_{2j+1}(H'_{n}-X)\leq \frac{4k+6}{8k+3}|X|+\frac{2k+3}{8k+3}$.
\end{lem}
%%%%%%%%%%%%%%%%%%%%%%%%%%%%%%%%%%%%%%%%%%%%%%%%%%%%%%%%%%%%%%%%%%%%%%%%%%%%%%%%%%%%%%%%%%%%%%%%%%%%%%%%%%%%%%%%%%%%%%%%
\proof
Let $X\subseteq V(H'_{n})$.

\begin{claim}%%%%%%%%%%%%%%%%%%%%%%%%%%%%%%%%%%%%%%%%%%%%%%%%%%%%%%%%%%%%%%%%%%%%%%%%%%%%%%%%%%%%%%%%%%%%%%%%%%%%%%%%%%%
\label{claim4.2.1.1}
For each $i~(1\leq i\leq 2n+1)$, $\sum _{0\leq j\leq k-1}c_{2j+1}(R_{i}-X)\leq \frac{4k+6}{8k+3}|V(R_{i})\cap X|+\frac{2k+3}{8k+3}$.
\end{claim}
%%%%%%%%%%%%%%%%%%%%%%%%%%%%%%%%%%%%%%%%%%%%%%%%%%%%%%%%%%%%%%%%%%%%%%%%%%%%%%%%%%%%%%%%%%%%%%%%%%%%%%%%%%%%%%%%%%%%%%%%
\proof
We first assume that $V(K_{i})\not\subseteq X$.
Then $R_{i}-X$ is connected.
Clearly we may assume that $\sum _{0\leq j\leq k-1}c_{2j+1}(R_{i}-X)=1$.
Then $|V(R_{i})\cap X|\geq 2$ because $|V(R_{i})|=2k+1$.
Hence $\sum _{0\leq j\leq k-1}c_{2j+1}(R_{i}-X)=1<\frac{4k+6}{8k+3}\cdot 2<\frac{4k+6}{8k+3}|V(R_{i})\cap X|+\frac{2k+3}{8k+3}$.
Thus we may assume that $V(K_{i})\subseteq X$.

Let $\alpha $ be the number of components of $R_{i}-V(K_{i})$ intersecting with $X$.
Since $\alpha \leq 2m+1$, we have $(8m+1)\alpha \leq (4m+2)(2m-1+\alpha )+2m+1$, and hence
$$
\alpha \leq \frac{4m+2}{8m+1}(2m-1+\alpha )+\frac{2m+1}{8m+1}=\frac{4k+6}{8k+3}(2m-1+\alpha )+\frac{2k+3}{8k+3}.
$$
Furthermore, $\sum _{0\leq j\leq k-1}c_{2j+1}(R_{i}-X)=c_{1}(R_{i}-X)\leq \alpha $ and $|V(R_{i})\cap X|=|V(K_{i})|+|(V(R_{i})-V(K_{i}))\cap X|\geq 2m-1+\alpha $.
Consequently we get $\sum _{0\leq j\leq k-1}c_{2j+1}(R_{i}-X)\leq \frac{4k+6}{8k+3}|V(R_{i})\cap X|+\frac{2k+3}{8k+3}$.
\qed

Assume for the moment that $V(R_{0})\not\subseteq X$.
Then $H'_{n}-X$ is connected.
Clearly we may assume that $\sum _{0\leq j\leq k-1}c_{2j+1}(H'_{n}-X)=1$.
Then $|X|\geq 2$ because $|V(H'_{n})|\geq 2k+1$.
Hence $\sum _{0\leq j\leq k-1}c_{2j+1}(H'_{n}-X)=1<\frac{4k+6}{8k+3}\cdot 2<\frac{4k+6}{8k+3}|X|+\frac{2k+3}{8k+3}$.
Thus we may assume that $V(R_{0})\subseteq X$.
Then clearly
\begin{align}
|&\C_{2j+1}(H'_{n}-X)|=\sum _{1\leq i\leq 2n+1}|\C_{2j+1}(R_{i}-X)|.\label{sec4-7}
\end{align}
By Claim~\ref{claim4.2.1.1} and (\ref{sec4-7}),
\begin{align*}
\sum _{0\leq j\leq k-1}c_{2j+1}(H'_{n}-X) &= \sum _{0\leq j\leq k-1}\left(\sum _{1\leq i\leq 2n+1}c_{2j+1}(R_{i}-X)\right)\\
&\leq \sum _{1\leq i\leq 2n+1}\left(\frac{4k+6}{8k+3}|V(R_{i})\cap X|+\frac{2k+3}{8k+3}\right)\\
&= \frac{4k+6}{8k+3}(|X|-|V(R_{0})|)+\frac{2k+3}{8k+3}(2n+1)\\
&= \frac{4k+6}{8k+3}(|X|-n)+\frac{2k+3}{8k+3}(2n+1)\\
&= \frac{4k+6}{8k+3}|X|+\frac{2k+3}{8k+3}.
\end{align*}
Thus we get the desired conclusion.
\qed

From Lemma~\ref{lem4.2.1}, we get the following proposition, which implies that if Conjecture~\ref{con1} is true, then the coefficient of $|X|$ in the conjecture is best possible.

\begin{prop}%%%%%%%%%%%%%%%%%%%%%%%%%%%%%%%%%%%%%%%%%%%%%%%%%%%%%%%%%%%%%%%%%%%%%%%%%%%%%%%%%%%%%%%%%%%%%%%%%%%%%%%%%%%%
\label{prop3.2}
For an integer $k\geq 3$ with $k\equiv 0~(\mbox{mod }3)$, there exist infinitely many graphs $G$ having no $\{P_{2},P_{2k+1}\}$-factor such that $\sum _{0\leq i\leq k-1}c_{2i+1}(G-X)\leq \frac{4k+6}{8k+3}|X|+\frac{2k+3}{8k+3}$ for all $X\subseteq V(G)$.
\end{prop}
%%%%%%%%%%%%%%%%%%%%%%%%%%%%%%%%%%%%%%%%%%%%%%%%%%%%%%%%%%%%%%%%%%%%%%%%%%%%%%%%%%%%%%%%%%%%%%%%%%%%%%%%%%%%%%%%%%%%%%%%


\begin{thebibliography}{99}
\bibitem{AAE}
J.~Akiyama, D.~Avis and H.~Era,
On a $\{1,2\}$-factor of a graph,
TRU Math. \textbf{16} (1980) 97--102.

\bibitem {D}
R.~Diestel,
``Graph Theory'' (4th edition), Graduate Texts in Mathematics \textbf{173},
Springer (2010).

\bibitem{K}
A.~Kaneko,
A necessary and sufficient condition for the existence of a path factor every component of which is a path of length at least two,
J. Combin. Theory Ser. B \textbf{88} (2003) 195--218.

\bibitem{KLS}
M.~Kano, C.~Lee and K.~Suzuki,
Path and cycle factors of cubic bipartite graphs,
Discuss. Math. Graph Theory \textbf{28} (2008) 551--556.

\bibitem{KMOO}
K.~Kawarabayashi, H.~Matsuda, Y.~Oda and K.~Ota,
Path factors in cubic graphs,
J. Graph Theory \textbf{39} (2002) 188--193.

\bibitem{LP}
M.~Loebl and S.~Poljak,
Efficient subgraph packing,
J. Combin. Theory Ser. B \textbf{59} (1993) 106--121.

\end{thebibliography}
\end{document}